\documentclass[twocolumn]{autart}

\usepackage{amsmath}
\usepackage{amssymb}
\usepackage{graphicx}
\usepackage{float}
\PassOptionsToPackage{hyphens}{url}\usepackage{url}
\usepackage{algorithm}
\usepackage{algpseudocode}

\usepackage{stackengine}

\usepackage{subcaption}

\setstackEOL{\\}
\makeatletter
\def\BState{\State\hskip-\ALG@thistlm}
\makeatother

\makeatletter
\def\endfigure{\end@float}
\def\endtable{\end@float}
\makeatother

\usepackage{stfloats}

\newtheorem{assumption}{Assumption}
\newtheorem{definition}{Definition}

\newcommand{\mbb}{\mathbb}
\renewcommand{\Re}{\mbb{R}}

\newcommand{\set}[2]{\left\{ #1\ \left| \ #2 \right. \right\}}

\newcommand{\innerprod}[2]{\langle{#1},{#2}\rangle}

\usepackage{bbm}
\usepackage{enumerate}
\usepackage{flushend}
\usepackage{setspace}
\allowdisplaybreaks

\begin{document}
\begin{frontmatter}

\title{Gradient-Bounded Dynamic Programming for Submodular and Concave Extensible Value Functions with Probabilistic Performance Guarantees\thanksref{footnoteinfo}}

\thanks[footnoteinfo]{Research is supported by SIA Food Union Management. A preliminary version of the results presented in this paper can be found in \cite{LEBEDEVETAL2020A}. These results have been extended in multiple directions: We provide a novel validation procedure, a probabilistic analysis and a more detailed case study.}

\author[First]{Denis Lebedev\thanksref{corrAuth}},
\author[First]{Paul Goulart},
\author[First]{Kostas Margellos}
\thanks[corrAuth]{Corresponding author.}

\address[First]{Department of Engineering Science, University of Oxford, Oxford, OX1 3PJ, United Kingdom (e-mail: \{denis.lebedev, paul.goulart, kostas.margellos\}@eng.ox.ac.uk)}


\begin{abstract}                
We consider stochastic dynamic programming problems with high-dimensional, discrete state-spaces and finite, discrete-time horizons that prohibit direct computation of the value function from a given Bellman equation for all states and time steps due to the ``curse of dimensionality''. For the case where the value function of the dynamic program is concave extensible and submodular in its state-space, we present a new algorithm that computes deterministic upper and stochastic lower bounds of the value function in the realm of dual dynamic programming. We show that the proposed algorithm terminates after a finite number of iterations. Furthermore, we derive probabilistic guarantees on the value accumulated under the associated policy for a single realisation of the dynamic program and for the expectation of this value. Finally, we demonstrate the efficacy of our approach on a high-dimensional numerical example from delivery slot pricing in attended home delivery.
\end{abstract}
\end{frontmatter}
\section{Introduction}
Dynamic programming (DP) is an established tool to solve many optimal control problems. Since the exact DP approach typically scales poorly with problem dimension, various approximation schemes have been developed. For linear multistage stochastic systems, for example, stochastic dual DP provides a remedy to the ``curse of dimensionality'' by constructing successively tighter upper and lower bounds to the exact value function of the DP \cite{PEREIRA1991}. However, current developments are mainly concerned with linear systems with piecewise-affine and piecewise-quadratic dynamics \cite{SHAPIRO2011,WARRINGTONETAL2017}. To the best of our knowledge, research on problems with discrete states has attracted only limited interest to date \cite{ZOU2019}. 

In this paper, we extend the stochastic dual DP approach to value functions that are both submodular and concave extensible over a discrete domain. For any DP whose value function possesses these properties, we present a new algorithm that computes deterministic upper bounds and stochastic lower bounds to the exact value function in the spirit of stochastic dual dynamic programming theory. Furthermore, we derive probabilistic guarantees on the value accumulated under the associated policy for a single realisation of the dynamic program and for the expectation of this value.
	
Value functions with the required properties for our approximation algorithm arise, for example, in network revenue management problems in a variety of applications like transportation, hospitality and appointment scheduling \cite{ZHANG2009,MEISSNER2012,SAURE2012}. One particular example is the revenue management problem in attended home delivery, i.e.\ the problem of finding the optimal pricing policy for a finite number of capacity-constrained delivery options. For this problem, \cite[Theorem 2]{LEBEDEVETAL2019B} showed that the value function is concave extensible and submodular under certain conditions. We provide a numerical example with problem parameters that stem from this application to show the effectiveness of our proposed algorithm. Our approach extends the approximate DP approach of \cite{YANG2017}, whose value function is approximated by an affine function, however without being accompanied by guarantees to generate a certain fraction of the optimal profit.
	
The paper is structured as follows: Section \ref{sec:problem} formulates our problem of interest and the assumptions that our work builds upon. In Section \ref{sec:approxalg}, we present a novel algorithm to compute approximately optimal policies for value functions over discrete state-spaces under assumptions on submodularity and concave extensibility. Section \ref{sec:approxres} derives deterministic upper bounds and stochastic lower bounds to the exact value function and shows convergence of the algorithm in a finite number of iterations. In Section \ref{sec:valalg}, we present an algorithm that validates the policy obtained in Section \ref{sec:approxalg} by computing sample profits, their empirical mean and their standard deviation. Section \ref{sec:valres} details our theoretical results on tail and expectation bounds of the sample profits obtained in Section \ref{sec:valalg}. In Section \ref{sec:gbdp_example}, we present a numerical example on a high-dimensional problem that is unsolvable by direct computation. Finally, we conclude in Section \ref{sec:conclusions} and provide directions for future research.

\emph{Notation: } For any $s\in \mathbb{N}$, let $1_s$ be a column vector of all zeros apart from the $s$-th entry, which equals 1. Furthermore, we define the convention that $1_0$ is a vector of zeros. Let $\mathbf{1}$ denote a vector of ones. Let $\innerprod{\cdot}{\cdot}$ denote the standard inner product of its arguments. let $\mathbb{E}$ denote the expectation operator, let $\Pr(\cdot)$ denote the probability of its argument and let $\mathbbm{1}(\cdot)$ denote the indicator function.

\section{Problem statement}\label{sec:problem}
We consider a discrete-space, discrete-time, finite horizon DP. Define discrete states $x \in X \subset \mathbb{Z}^n$ and continuous and/or discrete decision variables \mbox{$d\in D \subset \mathbb{Z}^a\times \Re^b$}. Define the set $S:=\{1,2,...,n\}$. Let the transition probability between two states $x$ and $y$ under decision $d$ be $P_{x,y}(d)$, where we require $P_{x,y}(d)\geq 0$ for all $(x,y,d)\in X\times X\times D$. For all $x \in X$, we impose that $\sum_{y \in Y_+(x)}P_{x,y}(d)=1$, where $Y_+(x):=\{x+1_s\}_{s\in S \cup \{0\}}$. This requirement implies that transitions in $x$ are only possible in the positive direction and by at most a unit step along one dimension. Such models are typical for order-taking processes 
\cite{YANG2017,YANG2016}. Furthermore, we define a finite time horizon $T:=\{1,2,\dots,\bar{t}\}$, a stage revenue function $g: \mathbb{Z}^{n}\times \mathbb{Z}^{n}\times (\mathbb{Z}^{a}\times\Re^b) \to \Re$ and a terminal cost function $C: \mathbb{Z}^n \to \Re$ to construct the following DP:
\begin{align}\label{eq:dp}
V_t(x) :=&\; \underset{d\in D}{\max}\left\{\sum_{y\in Y_+(x)}P_{x,y}(d)\left(g(x,y,d)+V_{t+1}(y)\right)\right\}\nonumber\\
&\;\forall (x,t)\in X\times T,\text{ where } \nonumber\\
V_{\bar{t}+1}(x):=&\;-C(x) \quad \forall x\in X.
\end{align}
It is not strictly necessary for $g$ to be independent of $t$ as long as the assumptions stated below can be satisfied. However, as our interest lies in time-independent problems and to ease notation, we ignore time-dependency of $g$ in this paper. To represent the DP more compactly, we define the Bellman operator $\mathcal{T}$ through the relationship
\begin{equation}\label{eq:dpcompact}
V_t = \mathcal{T}V_{t+1}, \text{ for all } t\in T.
\end{equation}
We next introduce several definitions needed to state the assumptions that we impose on the DP in \eqref{eq:dp}. 
\begin{definition}\label{de:submod}
A function $f:\mathbb{Z}^n\to \Re$ is \emph{submodular} if it satisfies
\begin{equation}\label{eq:submodular}
f(\max(y,z))+f(\min(y,z)) \leq f(y)+f(z)
\end{equation}
for all $(y,z) \in \mathbb{Z}^n \times \mathbb{Z}^n$, where the maximum and minimum are taken elementwise.
\end{definition}
 The following two definitions are commonly used in discrete convex analysis:
\begin{definition}\label{de:conc_clo}
Let $a \in \Re^n$ and $b \in \Re$. Then the \emph{concave closure} $\tilde{f}: \Re^n \to \Re \cup -\infty$ of a function \mbox{$f:\mathbb{Z}^n \to \Re \cup -\infty$} is defined as \cite[equation (2.1)]{MUROTA2001}
\begin{equation}\label{eq:concclo}
\tilde{f}(x):=\underset{a,b}{\inf}\set{\innerprod{a}{x} + b}{\innerprod{a}{y} + b \geq f(y) \hspace{4mm} \forall y \in \mathbb{Z}^n}.\nonumber
\end{equation}
\end{definition}
\begin{definition}\label{de:conc_ext}
A function $f: \mathbb{Z}^n \to \Re \cup -\infty$ is \emph{concave extensible} if and only if the evaluations of $f$ coincide with the evaluations of its concave closure $\tilde{f}$ \cite[Lemma 2.3]{MUROTA2001}, i.e.\ $f(x)=\tilde{f}(x), \text{ for all } x\in \mathrm{dom}(\tilde{f}\,)$.
\end{definition}
These definitions allow us to state the assumptions that we impose on the DP in \eqref{eq:dp}:
\begin{assumption}\label{as:terminal}
The function $C$ is submodular and concave extensible in $x$.
\end{assumption}
\begin{assumption}\label{as:preserve}
The Bellman operator preserves concave extensibility and submodularity of any concave extensible and submodular value function, i.e. if\, $V_{t+1}$ is submodular and concave extensible, then $V_t=\mathcal{T}V_{t+1}$ also has these properties for all $t\in T$.
\end{assumption}
In \cite[Theorem 2]{LEBEDEVETAL2019B}, it is shown that, under mild technical assumptions on the customer arrival rate, these assumptions are satisfied for the revenue management problem considered in Section \ref{sec:gbdp_example}.
	
\section{Value function approximation algorithm} \label{sec:approxalg}
We first state our proposed approximation procedure in Algorithm \ref{alg:GBDP} below and subsequently detail the individual algorithm steps. Inspired by stochastic dual DP techniques \cite{SHAPIRO2011}, the main idea of our algorithm is to alternate between generating sample paths in ``forward sweeps'' and refining the value function in ``backward sweeps''. The following sections describe this procedure in detail.

\begin{algorithm}[H]
\setstretch{1.05}
\caption{Proposed approximation algorithm}\label{alg:GBDP}
\begin{algorithmic}[1]
\State Initialise parameters: $X, D, P_{x,y}, T, g, C$ and $i_{\max}$
\State Initialise $Q_t^{0}(x) \gets \infty$, for all $(x,t)\in X\times T$
\State Initialise $Q_{\bar{t}+1}^0(x) \gets -C(x)$, for all $x \in X$
\For{$i\in \{1,2,\dots, i_{\max}\}$}
	\State $x_1^i \gets 0$
	\For{$t\in T$} \Comment ``Forward sweep''
		\State \Longunderstack[l]{$d_t^i \gets d^* \in \underset{d \in D}{\text{argmax}}\left\{\sum_{x_{t+1}^i\in Y_+(x_{t}^i)}P_{x_t^i,x_{t+1}^i}(d) \right.$\\ \\$\left.\hphantom{d_t^i \gets}\vphantom{\sum_{x_t}^1Q^1}\times\left(g(x_t^i,x_{t+1}^i,d)+Q_{t+1}^{i-1}(x_{t+1}^i)\right)\right\}$}
		\State $x_{t+1}^i \gets x_{t}^i + \underset{x_{t+1}^i}{\text{sample}}\left\{P_{x_t^i,x_{t+1}^i}\left(d_t^i\right)\right\}$
	\EndFor
	\State $l(i)\gets \sum_{t=1}^{\bar{t}} g(x_t^i,x_{t+1}^i,d_t^{i}) - C(x_{\bar{t}+1})$
	\While{$t > 1$} \Comment ``Backward sweep''
		\State $Z(x_{t+1}^i)\gets\{x_{t+1}^i+1_s+1_s'\}_{\begin{subarray}{l}
		 s\in S\cup\{0\},\\
		s'\in S\cup\{0\}\end{subarray}}$
		\If{$Q_{t+1}^{i-1}$ is submodular on $Z(x_{t+1}^i)$}
			\State $H^* \gets$ unique hyperplane through \hphantom{www} $\hphantom{wwwwww}\vphantom{sum^1}\left\{\left(y,(\mathcal{T}Q_{t+1}^{i-1})(y)\right)\right\}_{y \in Y_{+}(x_{t+1}^i)}$
		\Else
			\State $j^* \in \underset{j\in J_{t+1}^{i-1}}{\text{argmin}}\left\{\left(\mathcal{T}H_{t+1}^{j-1}\right)\left(x_{t+1}^i\right)\right\}$	
			\State $H^* \gets \mathcal{T}H_{t+1}^{j^*-1}$
  		\EndIf
  		\State $Q_t^{i}\gets \min\left\{H^*,Q_t^{i-1}\right\}$
  		\State $t \gets t-1 $ 
	\EndWhile
	\State $u(i)\gets Q_1^{i}(0)$ 
\EndFor
\end{algorithmic}
\end{algorithm}
\subsection{Initialisation}
We first initialise all parameters of the DP in \eqref{eq:dp} (step 1). Denote the maximum number of iterations by $i_{\max} \in \mathbb{N}$ and let $I:=\{0,1,\dots,i_{\max}\}$. Let the value function approximation $Q_t^{i}$ for all $(i,t)\in I\times T$ be the pointwise minimum of a finite number of affine functions, i.e.
\begin{equation}\label{eq:qminh}
Q_t^{i}(x):=\underset{j \in \{0,1,\dots,i\}}{\min}H_{t}^j(x), \text{ for all } x\in X,
\end{equation}
where $H_{t}^j: X \mapsto \Re$ describes a hyperplane, i.e.\ 
\begin{equation}
H_{t}^j(x):=\innerprod{a_{t}^j}{x} + b_{t}^j, \text{ for all } x\in X,
\end{equation}
with $a_{t}^j \in \Re^{n}, b_{t}^j \in \Re$ for all $(t,j)\in T\times I$. We characterise the set of supporting hyperplanes at $x$ as

\begin{equation}
J_t^{i}(x):=\underset{j \in \{0,1,\dots,i\}}{\text{argmin}}\left\{\innerprod{a_{t}^j}{x} + b_{t}^j \right\}
\end{equation}
for all $(x,i,t)\in X\times I \times T$. We construct $Q_t^i$ as a successively tighter upper bound of $V_t$ (as $i$ increases), i.e.\ $V_t(x)\leq Q_t^{i}(x)\leq Q_t^{i-1}(x)$ for all $(x,i,t)\in X\times (I\setminus \{0\})\times T$. In the $i$-th ''backward sweep``, $H_t^i$ is added to $Q_t^{i-1}$ for all $t\in T$ to form $Q_t^i$. To initialise $Q_t^0$, one could simply set $Q_t^0$ to be a single affine function with $a_t^0=0$ and $b_t^0=\infty$, such that $Q_t^0$ is indeed an upper bound to $V_t$ for all $t\in T$ (step 2). We discuss the possibility of closer initialisations in the context of our example in Section \ref{sec:gbdp_example}. We also initialise $Q_{\bar{t}+1}^i(x) := V_{\bar{t}+1}(x) = - C(x)$ for all $(x,i)\in X\times I$, which is a tight upper bound by the construction of the DP in \eqref{eq:dp} (step 3).

\subsection{``Forward sweep''}
Fix any iteration $i\in I\setminus\{0\}$. In each ``forward sweep'', we solve an approximate version of the Bellman equation in \eqref{eq:dp} forward in time, i.e.\ by replacing $V_t$ with its approximation $Q_t^{i-1}$ (step 7). Hence, we compute suboptimal $d_t^i$ for all $t \in T$ and simulate state transitions by sampling from the transition probability distribution given the  approximately optimal decisions (step 8). This defines a sample path $x_{t}^i $ for all $t \in T\cup \{\bar{t}+1\}$. At the end of each ``forward sweep'', we compute a stochastic lower bound on the total expected profit $V_1(0)$, which we denote by $l(i)$ for all $i \in I\setminus\{0\}$ (step 10). We show that this is indeed a stochastic lower bound in Section \ref{sec:valres}.

\subsection{``Backward sweep''}
Fix any iteration $i\in I$. In each ``backward sweep'', we first check if $Q_{t+1}^{i-1}$ is submodular on $Z(x_{t+1}^i)$ by computing the sign of \eqref{eq:submodular} for all possible pairs of points $(y,y') \in Z(x_{t+1}^i) \times Z(x_{t+1}^i)$, such that $y\neq y'$ (step 12). If the inequality in \eqref{eq:submodular} holds for all these points, we locally compute the exact DP stage problem on the set $Y_+(x_{t+1}^i)$, i.e.\ $\{\mathcal{T}Q_{t+1}^{i-1}(y)\}_{y \in Y_+(x_{t+1}^i)}$, to construct the hyperplane through $\left\{\left(y,(\mathcal{T}Q_{t+1}^{i-1})(y)\right)\right\}_{y \in Y_+(x_{t+1}^i)}$ (step 14). Then, the resulting added hyperplane is an upper bound to $V_t(x)$ for all $x \in X$, as shown in Section \ref{sec:approxres}. 

In the opposite case, we need to compute a submodular upper bound on $Q_{t+1}^{i-1}$, which is readily given by the hyperplanes from which $Q_{t+1}^{i-1}$ is constructed. Therefore, we select the hyperplane $H_{t+1}^{j^*-1}$ that minimises the value at the evaluation point $x_{t}^{i}$, which therefore locally creates the tightest upper bound (step 16). It may be possible to construct other submodular upper bounds to non-submodular $Q_{t+1}^{i-1}$, however, steps 16 and 17 of Algorithm \ref{alg:GBDP} offer a simple implementation. Finally, we update the value function approximation as the pointwise minimum of the approximation from the previous iteration and the newly constructed hyperplane (step 19). We also compute an upper bound, $u(i)$ for all $i\in I\setminus\{0\}$, on the total expected profit $V_1(0)$ (step 22). We show that this is indeed an upper bound in Section \ref{sec:approxres}.

\section{Approximation algorithm properties}\label{sec:approxres}
In this section, we show our main theoretical results on bounds on the exact value function and convergence properties of Algorithm \ref{alg:GBDP}. Proofs not included in this section can be found in the Appendix.

\begin{prop}\label{pr:QgeqV}
Under Assumptions \ref{as:terminal} and \ref{as:preserve}, the approximate value function is an upper bound to the exact finite horizon value function, i.e.\ $Q_t^i(x) \geq V_t(x)$ for all $(x,i,t)\in X\times I \times T$.
\end{prop}

\begin{cor}\label{co:u}
Under Assumptions \ref{as:terminal} and \ref{as:preserve}, the value of $u(i)$ is an upper bound to the exact total expected profit, i.e.\ $u(i) \geq V_1(0)$ for all $i\in I\setminus\{0\}$.
\end{cor}
\begin{pf}
This result follows immediately from Proposition \ref{pr:QgeqV} and by observing that $u(i)=Q_1^i(0)$ for all $i \in I\setminus\{0\}$ from step 22 of Algorithm \ref{alg:GBDP}.
\end{pf}
\begin{prop}\label{pr:l}
The value of $l(i)$ is a stochastic lower bound to the expected total profit, i.e.\ $\mathbb{E}[l(i)]\leq V_1(0)$ for all $i\in I\setminus\{0\}$.
\end{prop}
\begin{pf}
For any $i\in I\setminus\{0\}$, the value of $l(i)$ is obtained from suboptimal decisions $d_t^i$ for all $t\in T$, due to the use of  $Q_{t+1}^{i-1}$ instead of the exact (yet unavailable) $V_{t+1}$ in step 7 of Algorithm \ref{alg:GBDP}. It follows that $d_t^i$ is not a maximiser of the exact DP in \eqref{eq:dp} which, by the principle of optimality, implies that the expected value accumulated under this suboptimal policy will not be greater than the value obtained under the optimal policy. Hence, $\mathbb{E}[l(i)]\leq V_1(0)$ for all $i\in I\setminus \{0\}$.
\end{pf}
The stochastic dual DP algorithm converges asymptotically in $i$ to the exact value function \cite{SHAPIRO2011}. We can strengthen this result for our algorithm by exploiting the fact that the set of states $X$ is finite. Hence, the proposed algorithm converges in a finite number of steps under the following additional assumption.
\begin{assumption}\label{as:resampling}
Step 8 of Algorithm \ref{alg:GBDP} is repeated if the exact value function at the state-time pair $(x,t)\in X\times T$ has already been computed, i.e.\ if\, $Q_{t}^{i-1}(x_{t+1}^i) = V_{t}(x_{t+1}^i)$. If this is the case, a state for which the exact value function has not yet been reached is selected at random, with positive probability for all states for which the value function has not been computed exactly yet.
\end{assumption}
\begin{prop}\label{pr:converge}
Under Assumptions \ref{as:terminal}, \ref{as:preserve} and \ref{as:resampling}, the gap $u(i)-\mathbb{E}[l(i)]$ for all $i\in I\setminus\{0\}$ converges to $0$ in at most $\bar{t}|X|$ training iterations of Algorithm \ref{alg:GBDP}.
\end{prop}

Note that it is likely to take an unacceptably large number of iterations for the algorithm to converge to the exact value function due to the large number of states $|X|$. Since the value function is computationally expensive to calculate for all states, we seek to generate closer approximations at points that are likely to be visited, i.e.\ points on the sample path, and to use this information to save on approximation accuracy for less likely samples.

Our ultimate objective is to solve problems with large state spaces ($|X|\approx 10^{20}$) and long time horizons ($|T|\approx10^4$). In such scenarios, the need to resample the state as detailed in Assumption \ref{as:resampling} becomes negligible, because the required number of iterations to reach convergence is much larger than the maximum acceptable number of iterations. Therefore, from a practical point of view, we do not resample to satisfy Assumption \ref{as:resampling}. In this case, the proposed algorithm only converges asymptotically to the exact value function instead of in a finite number of steps, just as in stochastic dual DP \cite{SHAPIRO2011}.

\section{Proposed validation algorithm}\label{sec:valalg}
As noted in the previous section, absolute convergence of the approximate value function to the exact value function cannot be achieved for industry-sized problems due to the ``curse of dimensionality''. The performance of the algorithm, i.e.\ how close the stochastic lower bound $l(i)$ is to the deterministic upper bound $u(i)$ for any $i\in I\setminus\{0\}$, can only be validated statistically to a certain level of probabilistic confidence. To this end, we will generate a set of validation samples as described in Algorithm \ref{alg:validate} and detailed further below.
\begin{algorithm}[h]
\caption{Proposed validation algorithm}\label{alg:validate}
\begin{algorithmic}[1]
\State Compute approximation: $Q_t^{i_{\max}}$, for all $t\in T\setminus\{0\}\cup\{\bar{t}+1\}$
\State Initialise number of validation samples $k_{\max}$
\For{$k\in K:=\{1,2,\dots,k_{\max}\}$}
	\State $x_1^k \gets 0$
	\For{$t\in T$} \Comment ``Forward validation sweep''
		\State \Longunderstack[l]{$d_t^k \gets d^* \in \underset{d \in D}{\text{argmax}}\left\{\sum_{x_{t+1}^k\in X}P_{x_t^k,x_{t+1}^k}(d) \right.$\\ \\$\left.\hphantom{d_t^k \gets}\vphantom{\sum_{x_t}^1Q^1}\times \left(g(x_t^k,x_{t+1}^k,d)+Q_{t+1}^{i_{\max}}(x_{t+1}^k)\right)\right\}$}
		\State $x_{t+1}^k \gets x_{t}^k + \underset{x_{t+1}^k}{\text{sample}}\left\{P_{x_t^k,x_{t+1}^k}\left(d_t^k\right)\right\}$
	\EndFor
	\State $l_{\mathrm{v}}(k)\gets \sum_{t=1}^{T} g(x_{t}^k,x_{t+1}^k,d_t^{k})-C(x_{\bar{t}+1})$
\EndFor
\State $\bar{l}_{\mathrm{v}} \gets k_{\max}^{-1}\sum_{k=1}^{k_{\max}} l_{\mathrm{v}}(k)$
\State $\sigma_{\mathrm{v}} \gets \sqrt{(k_{\max}-1)^{-1}\sum_{k=1}^{k_{\max}} \left(l_{\mathrm{v}}(k)-\bar{l}_{\mathrm{v}}\right)^2}$
\end{algorithmic}
\end{algorithm}

We first compute the approximation obtained in Algorithm $\ref{alg:GBDP}$ (step 1). We denote the maximum number of validation samples by $k_{\max}\in \mathbb{N}$ and let $K:=\{1,2,\dots,k_{\max}\}$ (step 2). We then compute $k_{\max}$ ``forward validation sweeps'', where in each of them we use our most refined estimate, $Q_{t+1}^{i_{\max}}$ as our approximate value function (steps 5--8). After each sweep $k\in K$, we compute the stochastic lower bound $l_{\mathrm{v}}(k)$ on the total expected profit, similarly to $l(i)$ for all $i\in I\setminus \{0\}$ in Algorithm \ref{alg:GBDP} (step 9). We then compute the sample mean profit $\bar{l}_{\mathrm{v}}$ and unbiased empirical standard deviation $\sigma_{\mathrm{v}}$ of the set of sampled lower bounds $\{l_{\mathrm{v}}(k)\}_{k\in K}$ (steps 11--12). As detailed in the next section, these quantities will be used to generate one-sided confidence intervals, quantifying the performance of the decision policy associated with the approximate value function $Q_{t+1}^{i_{\max}}$.

\section{Validation algorithm properties}\label{sec:valres}
In this section we state the main theoretical properties of our validation procedure. The proofs can be found in the Appendix. We use $\{l_{\mathrm{v}}(k)\}_{k\in K}$, $\bar{l}_{\mathrm{v}}$ and $\sigma_{\mathrm{v}}$ from Algorithm \ref{alg:validate} to derive two different measures for the performance guarantee. The first is a probabilistic bound on the tail of the distribution of a single lower bound sample, i.e.\ a value for $l(k_{\max}+1)$ that is reached or exceeded with $1-\alpha$ confidence for a user-defined $\alpha\in(0,1)$. As we will see later in Section \ref{sec:gbdp_example}, this bound is not necessarily indicative of the expectation of $\bar{l}_{\mathrm{v}}$, since even under the profit-maximising decision policy, some variance will persist in $l(k_{\max}+1)$ from the randomness of the state transitions. Therefore, we also derive a bound on the expectation of the empirical sample mean $\bar{l}_{\mathrm{v}}$ that holds with confidence $1-\alpha^{\mathbb{E}}$, where $\alpha^{\mathbb{E}}\in(0,1)$ can be chosen by the user.

\subsection{Tail bounds}
In this section, we present two tail bounds of the distribution of $l_{\mathrm{v}}(k_{\max}+1)$. Let $[l_{-},l_{+}]$ denote the (finite) support of the distribution of $l_{\mathrm{v}}(k)$ for any $k\in K\cup \{k_{\max}+1\}$ and let $F_K$ denote the empirical cumulative distribution function of $\{l_{\mathrm{v}}(k)\}_{k\in K}$, i.e.\ \mbox{$F_K(l):=k_{\max}^{-1}\sum_{k\in K}\mathbbm{1}(l_{\mathrm{v}}(k)\geq l)$}. We derive two tail bounds of the distribution of $l_{\mathrm{v}}(k_{\max}+1)$ with a given confidence level $(1-\alpha)\in (0,1)$, which is mildly restricted for the first bound due to the next assumption.
\begin{assumption}\label{as:alphac}
Assume that $\alpha > \theta_{\mathrm{C}}:=\Pr(\sigma_{\mathrm{v}}=0)$.
\end{assumption}
The value of $\theta_{\mathrm{C}}$ will often be negligibly small, since $l_{\mathrm{v}}(k)$ for all $k\in K$ is highly unlikely to take identical values due to the typically high-dimensional state space and long time horizon. We show this later in Section \ref{sec:compbounds}. 
\begin{prop}\label{pr:oneoffbounds}
For any significance level $\alpha\in (0,1)$, $\Pr(l_{\mathrm{v}}(k_{\max}+1)\geq l^*) \geq 1-\alpha$ for all $l^* \in \{l_{\mathrm{C}},l_{\mathrm{D}}\}$, where:

\emph{(i)} under Assumption \ref{as:alphac}, $l_{\mathrm{C}}$ is the empirical Cantelli bound given by

\begin{equation}\label{eq:lc}
l_{\mathrm{C}}:= \bar{l}_{\mathrm{v}}- \sigma_{\mathrm{v}} \sqrt{\frac{(1-\alpha)(k_{\max}-1)}{(\alpha -\theta_{\mathrm{C}}) k_{\max}}}\text{ and }
\end{equation}
\emph{(ii)} $l_{\mathrm{D}}$ is the Dvoretzky-Kiefer-Wolfowitz bound given by
\begin{equation}\label{eq:ldkw}
l_{\mathrm{D}} := \mathrm{sup} \set{l \in [l_-,l_+]}{F_K(l)\leq \alpha-\theta_{\mathrm{D}}-\sqrt{\frac{\ln(\frac{1}{\theta_{\mathrm{D}}})}{2k_{\max}}}},
\end{equation}
where $\theta_{\mathrm{D}}\in (0,\alpha)$ is a user-defined parameter.
\end{prop}

For $l_{\mathrm{D}}$, we find the $\theta_{\mathrm{D}}$, which maximises the value of the bound, from the so-called Lambert $W$ function.
\begin{definition}\label{de:lambertw}
Let the Lambert $W$ function be implicitly defined as $W_{i}: \Re \to \Re$, such that $W_{i}(x)\exp(W_{i}(x))=x$ for $i\in \{0,-1\}$, where $W_{0}(x) > -1$ is called the principal branch and $W_{-1}(x) \leq -1$ is called the lower branch.
\end{definition}
\begin{lem}\label{le:opt_ld}
For any $\alpha \in (0,1)$, the value of $l_{\mathrm{D}}$ is maximised at 
\begin{equation}
\theta_{\mathrm{D}} = \min\left\{\alpha,\sqrt{\exp\left(W_{-1}\left(\frac{-1}{4k_{\max}}\right)\right)}\right\}.
\end{equation}
\end{lem}
The bounds $l_{\mathrm{C}}$ and $l_{\mathrm{D}}$, are termed after Cantelli's inequality \cite{CANTELLI1928} and the Dvoretzky-Kiefer-Wolfowitz \cite{MASSART1990} inequality, respectively. These inequalities are critical for showing that the bounds are indeed reached or exceeded with probability $1-\alpha$. By Proposition \ref{pr:oneoffbounds}, we can always choose the tighter, i.e.\ greater, of the two bounds and we will see later in Section \ref{sec:gbdp_example} that the selection of $\alpha$ and $k_{\max}$ influences which bound is preferred.

\subsection{Expectation bounds}
Similarly to the tail bounds, we now state our theoretical results on two bounds on the expectation of $\bar{l}_{\mathrm{v}}$, denoted by $\mathbb{E} \bar{l}_v$. To this end, we need to state two more technical assumptions.
\begin{assumption}\label{as:lmin}
Fix any $k\in \mathbb{N}$. The minimum value that $l_{\mathrm{v}}(k)$ can take is $l_{-}\geq 0$.
\end{assumption}
Assumption \ref{as:lmin} is weak, since it can always be satisfied by a linear transformation of $l_{\mathrm{v}}(k)$, such that $l_-\geq 0$.
\begin{prop}\label{pr:ebounds}
Fix any significance level $\alpha^{\mathbb{E}}\in (0,1)$. Then $\Pr(\mathbb{E}\bar{l}_{\mathrm{v}}\geq l^*) \geq 1-\alpha^{\mathbb{E}}$, for all $l^* \in \{l_{\mathrm{B}}^{\mathbb{E}},l_{\mathrm{D}}^{\mathbb{E}}\}$, where:

\emph{(i)} $l_{\mathrm{B}}^{\mathbb{E}}$ is the empirical Bernstein bound given by
\begin{equation}\label{eq:lbe}
l_{\mathrm{B}}^{\mathbb{E}}:=\bar{l}_{\mathrm{v}} - \sqrt{\frac{2\sigma_{\mathrm{v}}\ln(\frac{2}{\alpha^{\mathbb{E}}})}{k_{\max}}}-\frac{7(l_{+}-l_{-})\ln(\frac{2}{\alpha^{\mathbb{E}}})}{3(k_{\max}-1)} \text{ and }
\end{equation}
\emph{(ii)} under Assumption \ref{as:lmin}, $l_{\mathrm{D}}^{\mathbb{E}}$ is the Dvoretzky-Kiefer-Wolfowitz expectation bound given by
\begin{equation}\label{eq:ldkwe}
l_{\mathrm{D}}^{\mathbb{E}}:=\int_{l=0}^{\infty}1-\min\left\{1,F_K(l) + \sqrt{\frac{\ln(\frac{1}{\alpha^{\mathbb{E}}})}{2k_{\max}}}\right\}\mathrm{d}l.
\end{equation}
\end{prop}
The bounds $l_{\mathrm{B}}^{\mathbb{E}}$ and $l_{\mathrm{D}}^{\mathbb{E}}$ are termed after the empirical Bernstein \cite{MAURER2009} and Dvoretzky-Kiefer-Wolfowitz \cite{MASSART1990} inequalities, respectively. The proof of Proposition \ref{pr:ebounds}(i) is given in \cite{MAURER2009}. It can be shown that $l_{\mathrm{D}}^{\mathbb{E}}$ is at least as tight as Hoeffding's concentration bound \cite{HOEFFDING1963}, given by
\begin{equation}
l_{\mathrm{H}}^{\mathbb{E}}:=\bar{l}_{\mathrm{v}}-(l_{+}-l_{-})\sqrt{\frac{\ln(1/\alpha^{\mathbb{E}})}{2k_{\max}}}.
\end{equation}
In fact, under an additional technical assumption, we show that the expectation Dvoretzky-Kiefer-Wolfowitz bound is stricly better than Hoeffding's bound.
\begin{assumption}\label{as:alphae}
We assume that $\alpha$ and $k_{\max}$ are chosen to satisfy \mbox{$\sqrt{\ln(1/\alpha^{\mathbb{E}})/(2k_{\max})}>k_{\max}^{-1}$}.
\end{assumption}
Assumption \ref{as:alphae} is very mild, since even for only a single observation $k_{\max}=1$, the critical value of $\alpha^{\mathbb{E}}$ would be $e^{-2}\approx13.5 \%$, which is much larger than typical significance levels, e.g.\ 5\% or 1\%. Taking any smaller value of $\alpha^{\mathbb{E}}$ than the critical value will ensure that Assumption \ref{as:alphae} is always satisfied. Furthermore, for $k_{\max} > 1$, the constraint on $\alpha^{\mathbb{E}}$ is even less restrictive.
\begin{prop}\label{pr:DKW_H}
Under Assumptions \ref{as:lmin} and \ref{as:alphae}, the Dvoretzky-Kiefer-Wolfowitz expectation bound is strictly tighter than Hoeffding's concentration bound, i.e.\ $l_{\mathrm{D}}^{\mathbb{E}} > l_{\mathrm{H}}^{\mathbb{E}}$ for all $\alpha^{\mathbb{E}}\in (0,1)$.
\end{prop}
Finally, we note that other bounds have also been proposed in the literature, e.g.\ \cite{SHAPIRO2011} assumes that the distribution of $\bar{l}_{\mathrm{v}}$ is Gaussian and determines confidence intervals based on the corresponding standard score, i.e.\ a Gaussian lower bound on the expectation of $\bar{l}_{\mathrm{v}}$ would be
\begin{equation}\label{eq:lge}
l_{\mathrm{G}}^{\mathbb{E}}:=\bar{l}_{\mathrm{v}} -z(\alpha^{\mathbb{E}})\frac{\sigma_{\mathrm{v}}}{\sqrt{k_{\max}}},
\end{equation}  where $z(\alpha^{\mathbb{E}})$ is the standard score of the Gaussian distribution (in fact, Student's t-distribution, especially for small sample sizes $k_{\max}$, since the true variance of the underlying distribution of $\bar{l}_{\mathrm{v}}$ is approximated by $\sigma_{\mathrm{v}}^2$). We compare this with our proposed bounds in Section \ref{sec:gbdp_example}.
\section{Numerical example}\label{sec:gbdp_example}
We demonstrate our algorithm on an example of the so-called revenue management problem in attended home delivery. The objective is to price delivery time windows, called ``slots'', dynamically over a finite time horizon to control the customer purchasing process to maximise profits while ensuring that all orders can still be fulfilled.

In this problem, $S$ is the set of delivery slots and the components of $x$ are the number of orders placed in every delivery slot. The feasible set of states $X$ is defined by the maximum state vector $\bar{x}$, i.e.\ $X:=\set{x\in \mathbb{Z}^n}{0\leq x \leq \bar{x}}$. The set of delivery slot price vectors is $D := \set{d\in \Re^n}{d_s \in [\underline{d}, \bar{d}], s=1,2,\dots,n}$. Customer choice follows a multinomial logit model \cite{DONGETAL2009}:
\begin{align}
P_{x,x}(d)&:=(1-\lambda)+ \frac{\lambda}{\sum_{k\in S}\exp(\beta_c+\beta_k+\beta_dd_{k})+1}\nonumber,\\[1ex]
P_{x,x+1_s}(d)&:= \frac{\lambda\exp(\beta_c+\beta_s+\beta_d d_{s})}{\sum_{k\in S}\exp(\beta_c+\beta_k+\beta_dd_{k})+1}
\end{align}
for all $(x,d,s)\in X\times D \times S$, where $\lambda\in (0,1)$ is the probability that a customer arrives on the booking website, $\beta_c \in \Re$ denotes a constant offset, $\beta_s \in \Re$ represents a measure of the popularity for all delivery slots $s\in S$ and $\beta_d < 0$ is a parameter for the price sensitivity.  More details on the estimation of these parameters can be found in \cite{YANG2016}. The average revenue of an order is $r$ and the length of the time horizon, representing the booking period, is $\bar{t}$. The cost function $C$ represents the delivery cost for all lists of orders $x\in X$ accumulated at the end of the booking period. The challenge is to price the slots dynamically to maximise profits, which corresponds to solving a DP of the form of \eqref{eq:dp}, where $g(x,y,d):=r+d_s$ if $y=x+1_s$ for all $s \in S$ and otherwise, $g(x,y,d):=0$, i.e.\ the stage revenue is the average revenue plus delivery price for slot $s$ if slot $s$ is chosen and otherwise, it is zero. The DP in our numerical example takes the parameters in Table \ref{tab:alg01par} below, adapted from a real-world, multi-subarea case study by \cite{YANG2017} to a single delivery subarea scenario. Furthermore, we also adopt the customer choice parameters $\left(\beta_c,\beta_d,\{\beta_s\}_{s\in S}\right)$ from that paper.
\begin{table}[H]
\caption{Numerical example parameters.}
\begin{center}
\renewcommand{\arraystretch}{1}
\begin{tabular}{r|l}\label{tab:alg01par}
	$S$ & $\{1,2,\dots,17\}$ \\
	$\bar{x}$ & $[6,6,\dots,6]$ \\
	$\lambda$ & $0.008$ \\
	$\left[\underline{d},\bar{d}\,\right]$ & $[\textrm{\pounds}0,\textrm{\pounds}10]$ \\
	$r$ & \pounds$34.53$ \\
	$\bar{t}$ & $6990$ \\
	$C(x)$ & $\textrm{\pounds}0.083\times \innerprod{\mathbf{1}}{x}$ if $x\in X$ and $\infty$ otherwise
\end{tabular}
\end{center}
\end{table}
We have chosen $C(0)=0$, i.e.\ we ignore fixed costs, which have no effect on the pricing policy. Notice that for direct value function computation we would have to evaluate the Bellman equation in \eqref{alg:GBDP} for all $(x,t) \in X\times T$, i.e.\ $(6+1)^{17}\times 6990 \approx 1.6\times 10^{18}$ evaluations in our example. This is prohibitively large for any available computer. Hence, we use an approximate algorithm.

For this type of DP, \cite[Theorem 2]{LEBEDEVETAL2019B} showed that the Bellman operator preserves strict submodularity, i.e.\ the condition in \eqref{eq:submodular} holds with strict inequality, if a small enough $\lambda > 0$ is chosen. However, in this problem the terminal condition $C$ is only weakly submodular in $x$. In fact, it is an affine function of $x$ and hence modular. We assume that $\lambda=0.008$ from Table \ref{tab:alg01par} is small enough to satisfy Assumption \ref{as:preserve} in this problem. We note however, that in the absence of strict submodularity, we compromise on the theoretical guarantee that $u(i)$ (Corollary \ref{co:u}) is an upper bound to the exact value function.

To speed up computation, we initialise $Q_{t}^0$ for all $t\in T$ using the fixed point of DP, $V^*$, which is a known upper bound to the exact value function at any $(x,t)\in X\times T$, i.e.\ $V^*(x)\geq V_t(x)$. This is always the case when $\mathcal{T}$ in \eqref{eq:dpcompact} is a monotone operator \cite[Chapter 3]{BERTSEKAS2012}. In \cite{LEBEDEVETAL2019A}, it is shown that the fixed point is given analytically as
\begin{equation}\label{eq:fixedpoint}
V^*(x):= (\bar{d}+r)\innerprod{\mathbf{1}}{\bar{x}-x}-C(\bar{x}),\text{ for all } x\in X.
\end{equation}
Hence, we use this result to set $Q_{t}^0(x)=V^*(x)$ instead of $\infty$ for all $(x,t)\in X\times T$. Note that the fixed point in \eqref{eq:fixedpoint} is an affine function, so the initialiser has low complexity, i.e.\ only one affine function describes $Q_{t}^0$.

\subsection{Computation of approximate value function}
We implement Algorithm \ref{alg:GBDP} in Julia \cite{BEZANSON2017} and run $i_{\max}=100$ iterations on an i7-8565U CPU at 1.80 GHz processor base frequency and with 16GB RAM. The run time on our machine is 25 mins, 48 sec. In each iteration $i\in \{1,2,\dots i_{\max}\}$, we compute the upper bound on the expected profit $u(i)$ (Corollary \ref{co:u}) and the stochastic lower bound $l(i)$ (Proposition \ref{pr:l}), corresponding to the sample profit obtained in a single ``forward sweep'' of Algorithm \ref{alg:GBDP}, as discussed in Section \ref{sec:approxalg}. We also compute the cumulative moving average of the sample profits, i.e.\ $i^{-1}\sum_{j=1}^i l(j)$, which tends to the expected value of the stochastic lower bound $l(i)$ (Proposition \ref{pr:l}) as $i$ increases. Fig.\ \ref{fig:converge} shows how these bounds develop over 100 iterations. We make the following observations:
\begin{figure}[h]
\centering
  \includegraphics[width=\linewidth,trim={0cm 0cm 0cm 0cm},clip]{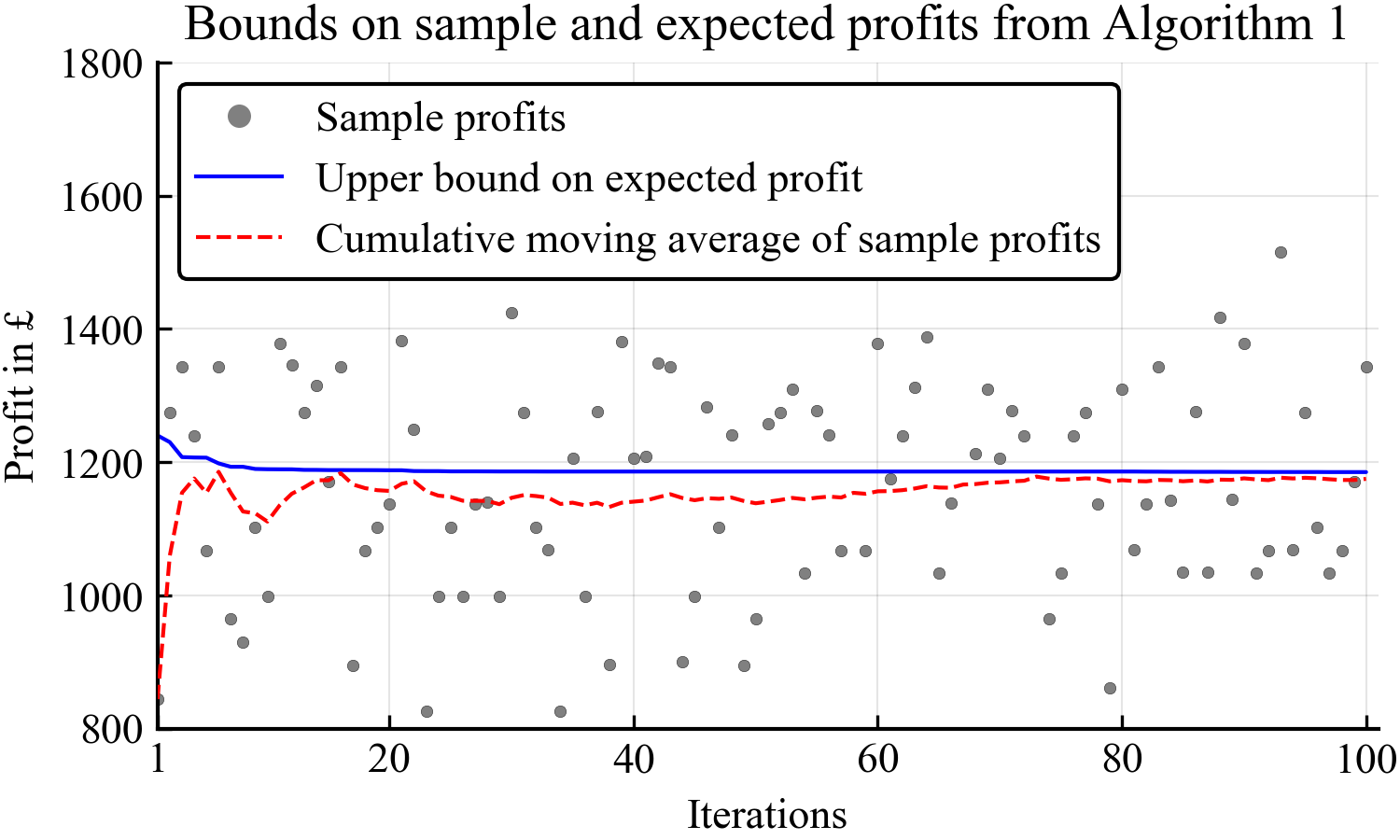}
  \caption{Plots of sample profits (grey dots), upper bound $u(i)$ (blue solid line, Corollary \ref{co:u}) and stochastic lower bound $l(i)$ (red dashed line, Proposition \ref{pr:l}) on expected profit.}
  \label{fig:converge}
\end{figure}

\emph{(1)} The upper bound converges within 10 iterations.

\emph{(2)} The cumulative moving average of the sample profits converges after about 100 iterations. We conjecture that this is due to the fact that Algorithm \ref{alg:GBDP} refines the value function approximation iteratively for all time steps. However, not all refinements propagate through all time steps of the DP. Therefore, each iteration step has a direct influence on the stochastic lower bound, which depends on the value function approximation at all time steps, while the upper bound might be unchanged as it only depends on the approximation at the first time step.

\emph{(3)} The sample profits have high variance at all iterations since the customer choice process is random.

We investigate the influence of additional iterations, by comparing the performance of the pricing policies obtained after 1 and after 100 iterations of Algorithm \ref{alg:GBDP}. In particular, for both iteration counts, we simulate 1,000 booking periods (by using the ``forward sweep'' of Algorithm \ref{alg:GBDP}) and we compute the sample profits obtained in each period. The resulting histogram of sample profits is shown in Fig.\ \ref{fig:histogram} below, where we observe the following:
\begin{figure}[b]
\centering
\begin{subfigure}{\linewidth}
\centering
\includegraphics[width=\textwidth]{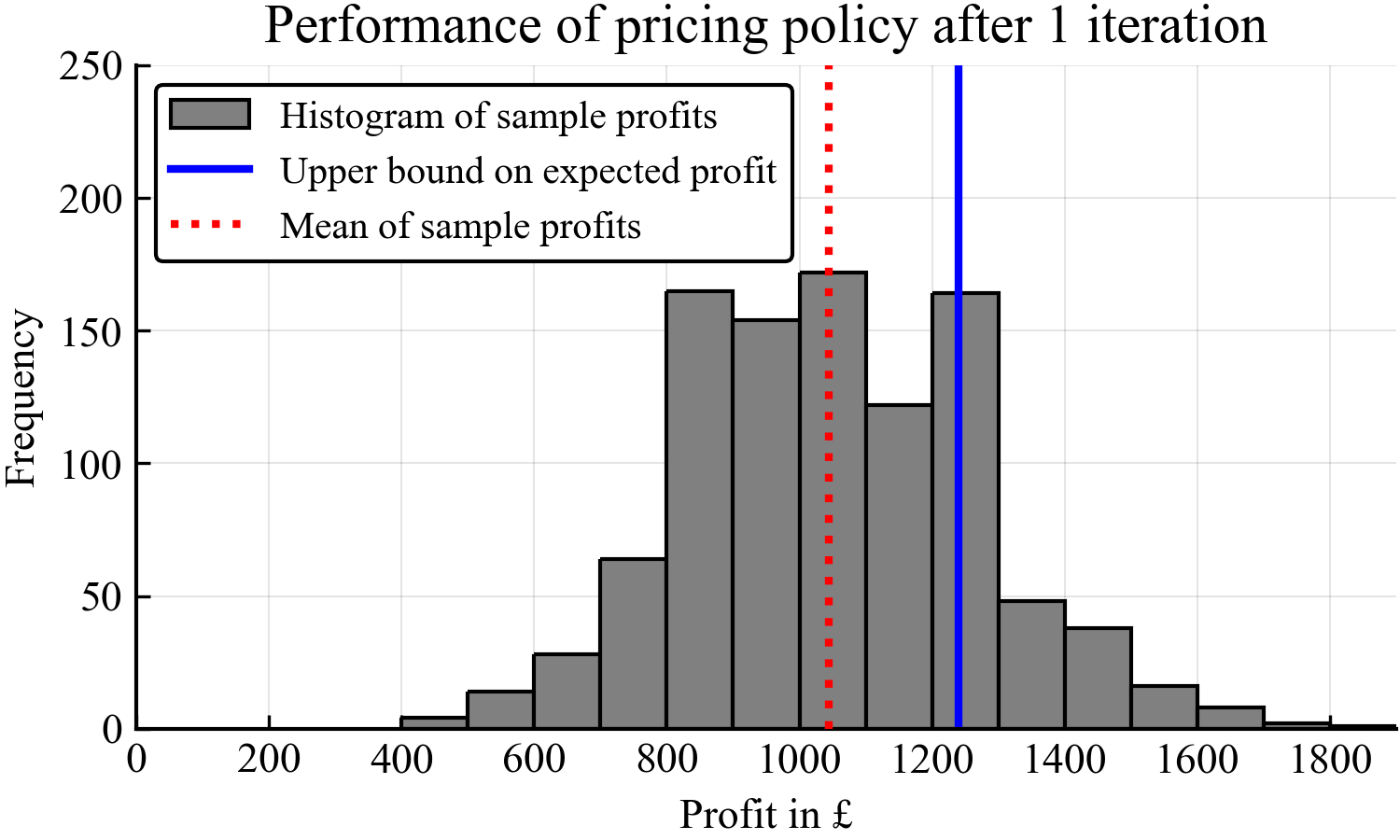}
\caption{After 1 iteration.}
\end{subfigure}\\
\hspace{0cm}
\begin{subfigure}{\linewidth}
\centering
\includegraphics[width=\textwidth]{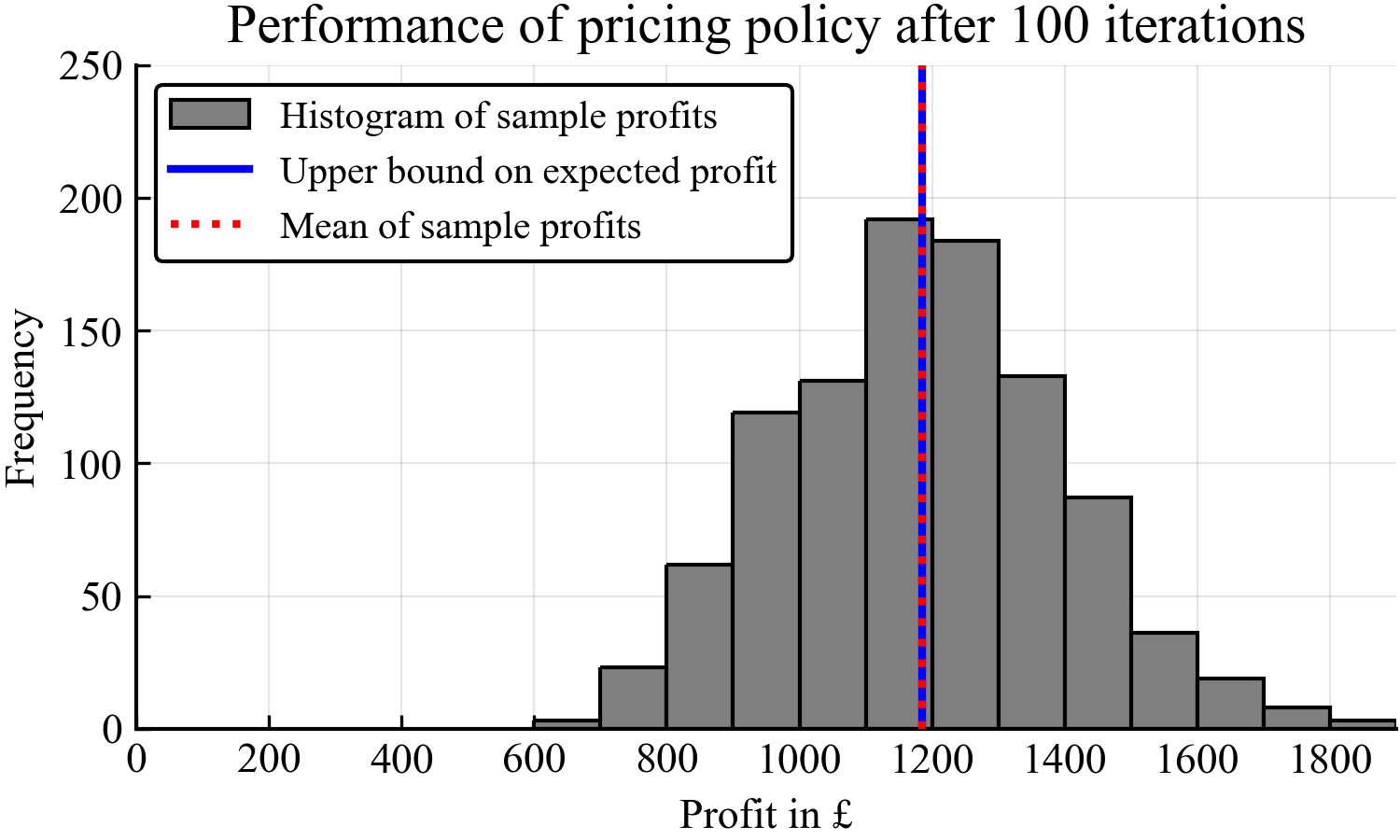}
\caption{After 100 iterations.}
\end{subfigure}
\caption{Histograms of 1,000 sample profits after 1 (a) and 100 (b) iterations of Algorithm \ref{alg:GBDP}.}
\label{fig:histogram}
\end{figure}

\emph{(1)} The mean sample profit increases from \pounds 1,044 after 1 iteration by 13.5\% to \pounds 1,185 after 100 iterations.

\emph{(2)} The gap between upper bound on expected profit, $u(i)$, and the empirical mean of 1,000 samples of $l(i)$, decreases from \pounds 200 to approximately \pounds 0. 

\emph{(3)} The histograms show that the sample profits are more concentrated around their empirical mean after 100 iterations than after 1 iteration. This indicates that the variance of the sample profits can be decreased by increasing the number of algorithm iterations.

We now validate our results by obtaining probabilistic performance guarantees for the approximate policy obtained after 100 iterations of Algorithm \ref{alg:GBDP}. In particular, we compute confidence intervals on the stochastic lower bound using validation samples $\{l_{\mathrm{v}}(k)\}$ for all $k\in K$ and the two pairs of bounds defined in Algorithm \ref{alg:validate}.

\subsection{Computation of bounds}\label{sec:compbounds}
We first compute the tail bounds on the value of profit obtained by a single sample under the approximate policy after 100 iterations of Algorithm \ref{alg:GBDP}. To this end, we compute the empirical Cantelli bound $l_{\mathrm{C}}$ from \eqref{eq:lc} and the Dvoretzky-Kiefer-Wolfowitz tail bound $l_{\mathrm{D}}$ from \eqref{eq:ldkw}. Due to the 17-dimensional state-space and long time horizon of 6990 steps, $\theta_{\mathrm{C}}\approx0$. To see this, upper bound the probability of the most likely event at every stage, namely no order being placed by
$P_{x,x}(d)\leq P_{x,x}(\mathbf{1}\underline{d}) \approx 0.9951$. Due to time-independence of the transition probabilities, we can exponentiate this number by the number of time steps in the DP to obtain the probability of 0 orders at the end of the booking period. This needs to happen for all $k_{\max}$ (independent) validation samples, hence we again exponentiate this number by $k_{\max}$. This gives us the probability of all validation samples having 0 orders. This is the most likely, but only one of $|X|$ states, so we multiply this number by $|X|$ to obtain $\Pr(\sigma_{\mathrm{v}}=0)\leq |X|P_{x,x}(\mathbf{1}\underline{d})^{\bar{t}k_{\max}}\approx 11^{17} \times 1.748\times 10^{-14000}\approx 0$.

\begin{figure*}[t]
\centering
\begin{subfigure}{0.32\linewidth}
\centering
\includegraphics[width=\textwidth]{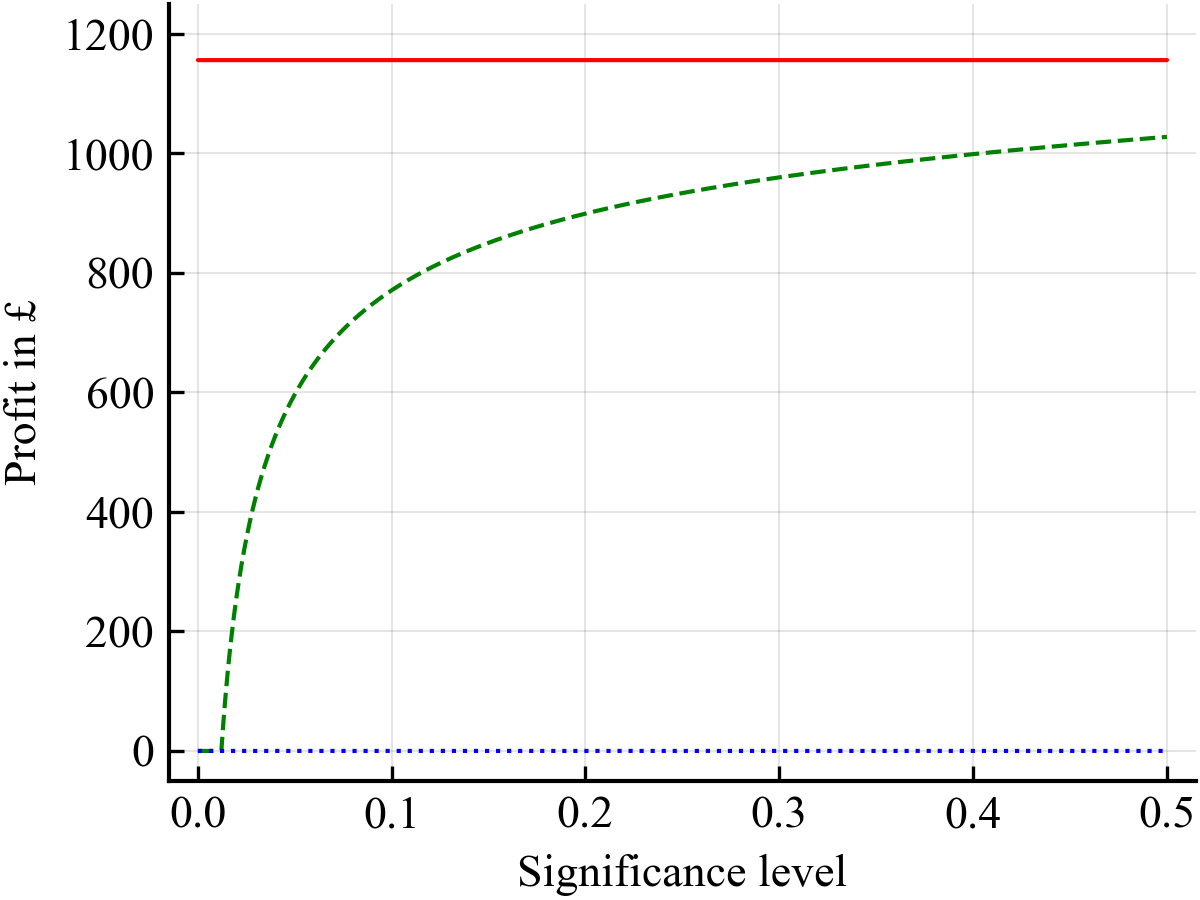}
\caption{10 samples.}
\end{subfigure}
\hspace{0cm}
\begin{subfigure}{0.32\linewidth}
\centering
\includegraphics[width=\textwidth]{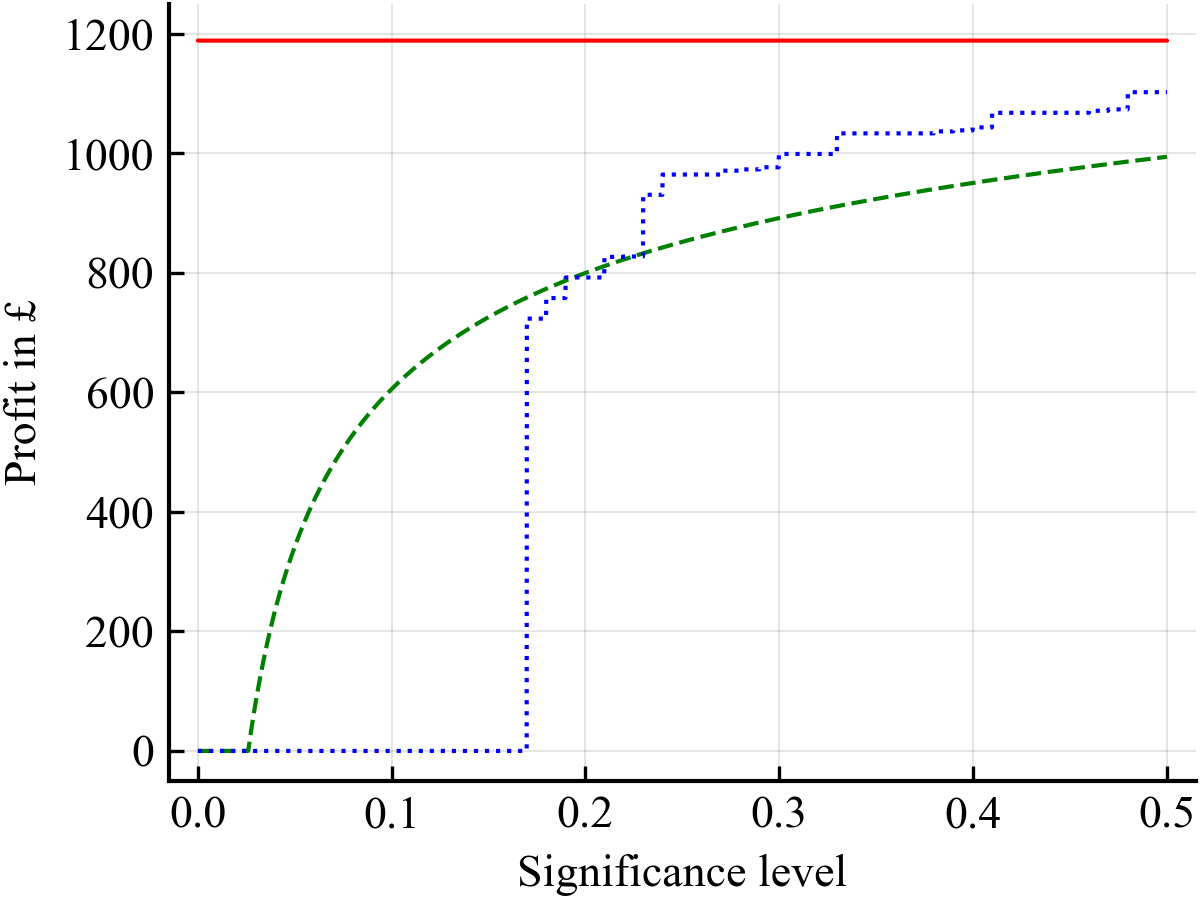}
\caption{100 samples.}
\end{subfigure}
\hspace{0cm}
\begin{subfigure}{0.32\linewidth}
\centering
\includegraphics[width=\textwidth]{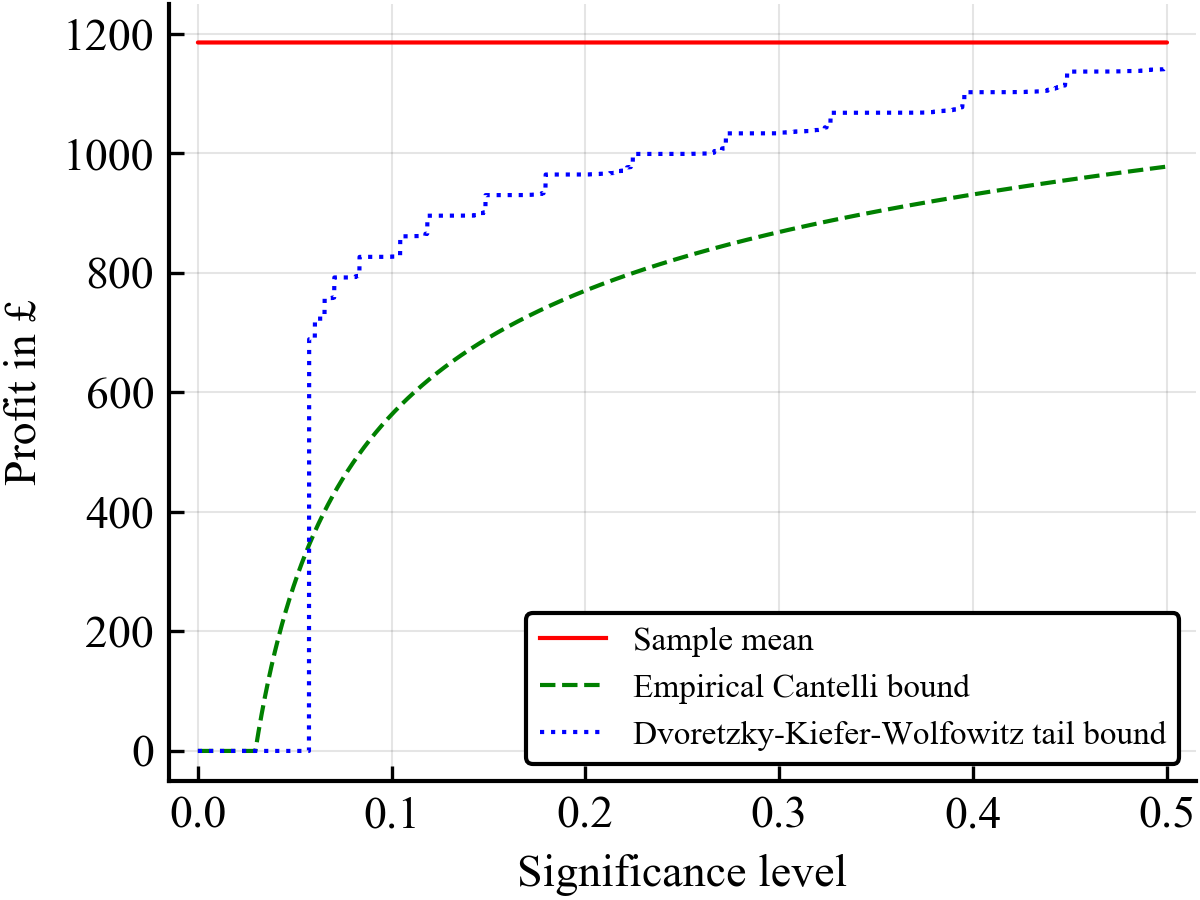}
\caption{1000 samples.}
\end{subfigure}
\caption{Probabilistic bounds -- Empirical Cantelli bound (green dashed line, Proposition \ref{pr:oneoffbounds}(i)) and Dvoretzky-Kiefer-Wolfowitz tail bound (blue dotted line, Proposition \ref{pr:oneoffbounds}(ii)) -- on the profit obtained from a single validation sample as functions of various significance levels for 10 (a), 100 (b) and 1000 (c) validation samples.}
\label{fig:oneoff}
\end{figure*}

As we see in Fig.\ \ref{fig:oneoff}, the tail bounds do not converge to the sample average, since there is an inherent variance in the customer choice model. This can be seen by inspecting the high variation of the sample profits in Fig.\ \ref{fig:converge} for all iteration steps. In Fig.\ \ref{fig:oneoff}, notice that both bounds create similar guarantees for different sample sizes: For $k_{\max}=10$ iterations and with probability 90\% the empirical Cantelli bound yields approximately \pounds800. A similar profit is obtained for $k_{\max}=1,000$ iterations and with probability 90\% using the Dvoretzky-Kiefer-Wolfowitz tail bound with optimal parameter $\theta_{\mathrm{D}}\approx \min\{\alpha,4.8\times 10^{-3}\}$ from Lemma \ref{le:opt_ld}.
\begin{figure*}[t]
\centering
\begin{subfigure}{0.32\linewidth}
\centering
\includegraphics[width=\textwidth]{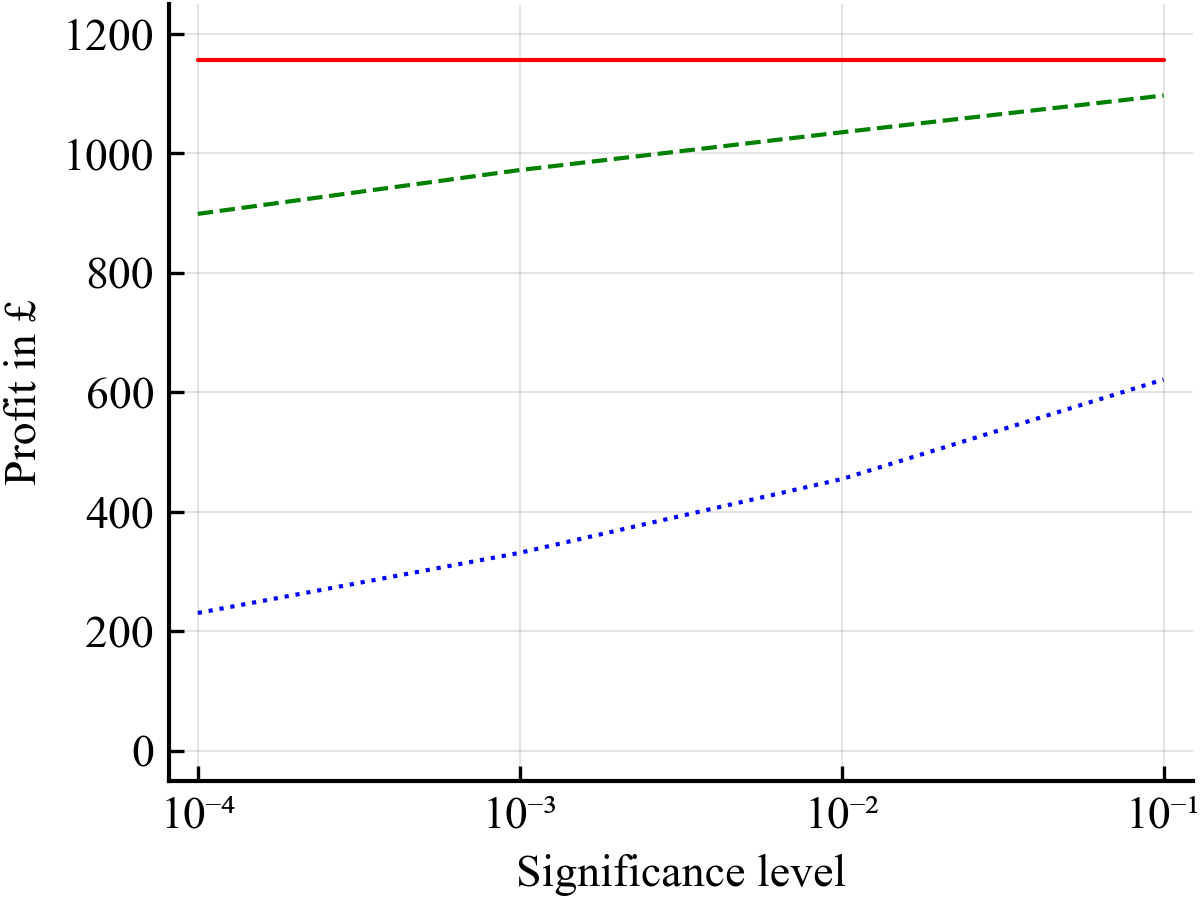}
\caption{10 samples.}
\end{subfigure}
\hspace{0cm}
\begin{subfigure}{0.32\linewidth}
\centering
\includegraphics[width=\textwidth]{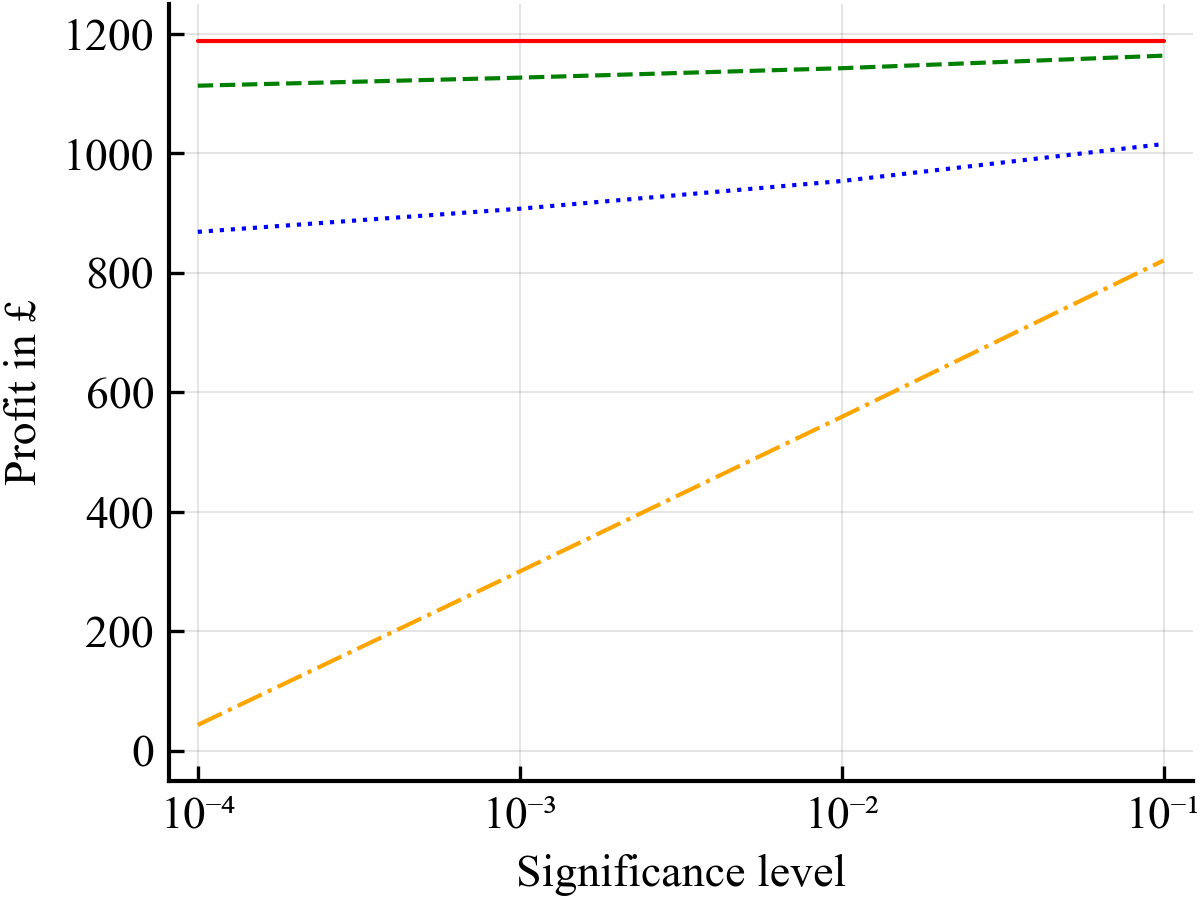}
\caption{100 samples.}
\end{subfigure}
\hspace{0cm}
\begin{subfigure}{0.32\linewidth}
\centering
\includegraphics[width=\textwidth]{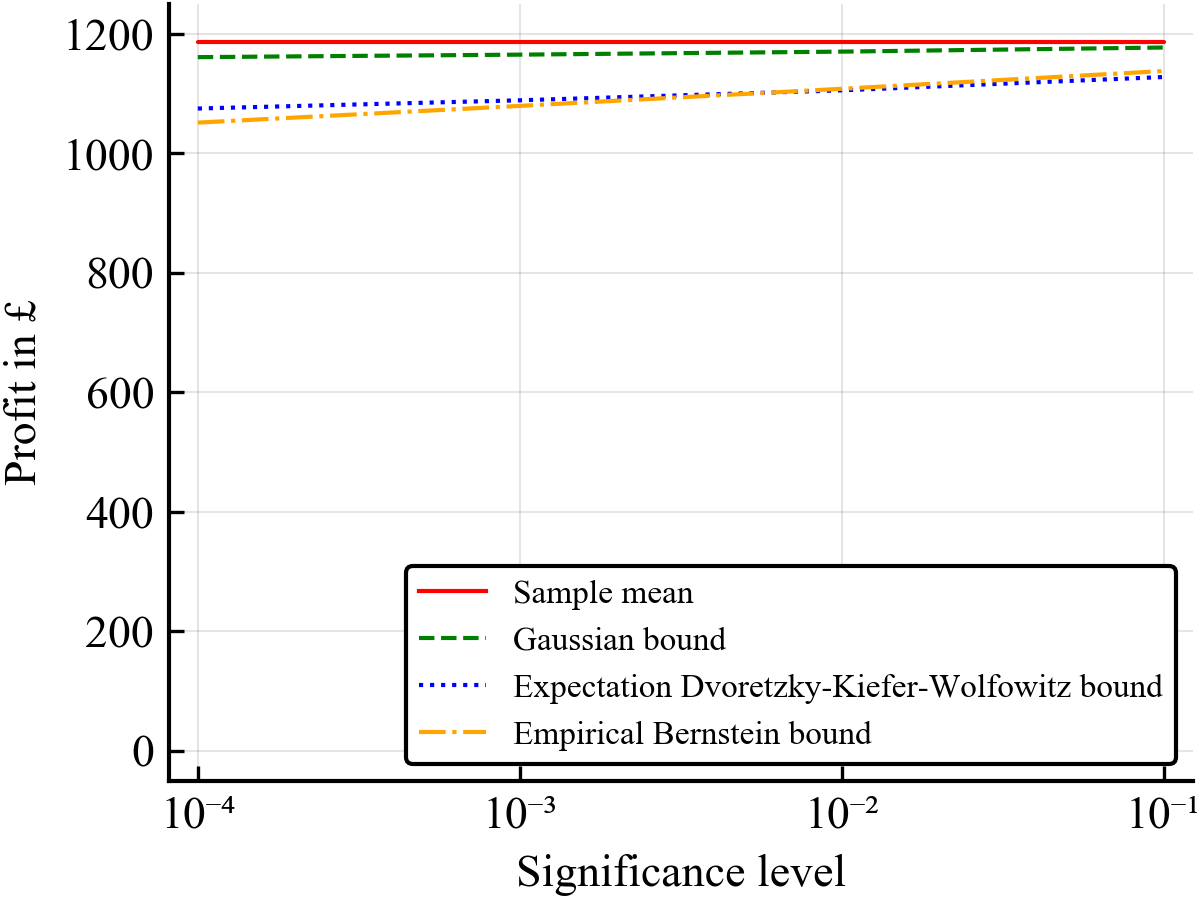}
\caption{1000 samples.}
\end{subfigure}
\caption{Probabilistic bounds -- Gaussian bound (green dashed line, \eqref{eq:lge}), Expectation Dvoretzky-Kiefer-Wolfowitz bound (blue dotted line, Proposition \ref{pr:ebounds}(ii)) and Empirical Bernstein bound (orange dash-dotted line, Proposition \ref{pr:ebounds}(i)) -- on the expectation of the profit obtained as functions of various significance levels and for 10 (a), 100 (b) and 1000 (c) validation samples.}
\label{fig:exp}
\end{figure*}

To get more meaningful measures of the convergence of Algorithm \ref{alg:GBDP}, we now compute the bounds on the expectation of the profit obtained after 100 iterations of Algorithm \ref{alg:GBDP}. To this end, we compute the empirical Bernstein bound $l_{\mathrm{B}}^{\mathbb{E}}$ from \eqref{eq:lbe} and the expectation Dvoretzky-Kiefer-Wolfowitz bound $l_{\mathrm{D}}^{\mathbb{E}}$ from \eqref{eq:ldkwe}. We also compute the Gaussian bound $l_{\mathrm{G}}^{\mathbb{E}}$ from \eqref{eq:lge}, following \cite{SHAPIRO2011}. This final bound relies on the additional assumption that the exact distribution of the mean sample profit $\bar{l}_{\mathrm{v}}$ is Gaussian. This is only asymptotically true by a Central Limit Theorem (see \cite[Proposition 2.16]{HAJEK2015}). As seen in Fig.\ \ref{fig:exp}, the Gaussian bound is always the most optimistic, however not accompanied by theoretical guarantees. Moreover, for large validation samples sizes $k_{\max}$, the gap with the empirical Bernstein bound and with the expectation Dvoretzky-Kiefer-Wolfowitz bound is small. 

Therefore, we suggest to sacrifice a small percentage of the profit gained by the Gaussian bound in return for the probabilistic guarantees of the other two bounds. Out of the  empirical Bernstein bound and the expectation Dvoretzky-Kiefer-Wolfowitz bound, the earlier only tends to perform better for large sample sizes ($k_{\max}=1000$) and at the same time not overly restrictive significance levels $(\alpha^{\mathbb{E}}>10^{-2})$. In all other cases, the expectation Dvoretzky-Kiefer-Wolfowitz bound is preferred. Since one typically wants to tighten significance levels as sample sizes increase, we anticipate that the expectation Dvoretzky-Kiefer-Wolfowitz bound will dominate the empirical Bernstein bound in many practical situations. Note that we omit the empirical Bernstein bound in Fig.\ \ref{fig:exp}(a) since its negative values are not meaningful.

\section{Conclusions and future work}\label{sec:conclusions}
We presented a new algorithm to compute approximate solutions to dynamic programs with submodular, concave extensible value functions. We derived deterministic upper and stochastic lower bounds to the exact value function. We showed that our algorithm converged in a finite number of iterations and we derived tail and expectation bounds for the stochastic lower bound. Finally, we demonstrated our results in an example of the revenue management problem in attended home delivery.
	
To the best of our knowledge, the proposed algorithm is the first to provide an upper bound on the expected profit for this class of problems. Comparing upper with stochastic lower bound, we can quantify the profit generation efficiency of our algorithm. This is a benchmark for other algorithms, possibly with weaker theoretical guarantees, but better practical performance. 

Finally, the gap between upper and stochastic lower bound will allow to quantify and optimise the trade-off between quality of approximation and computational cost \textit{dynamically} as an application runs. For example, in the revenue management problem in attended home delivery, the terminal condition of the dynamic program is an approximation to the intractable capacitated vehicle routing problem with time windows \cite{TOTH2014}. As time in the booking horizon progresses, orders come in and reveal the location of customers. This information could be used to update the terminal condition of the DP and hence, the pricing policy by re-running our algorithm.

\begin{ack}
We gratefully acknowledge the helpful discussions with Michael Garstka, Department of Engineering Science, University of Oxford, on the Julia implementation of our algorithm.
\end{ack}
\bibliographystyle{plain}
\bibliography{Algorithm02_Ref03}
\appendix
\section{Appendix}\label{ap:proofs}
\subsection{Proof of Proposition \ref{pr:QgeqV}}
We show this result by induction on $t$. In the base case (the terminal condition), $Q_{\bar{t}+1}^i(x) := V_{\bar{t}+1}(x) = - C(x)$ for all $(x,i)\in X\times I$, which satisfies the proposition trivially by Assumption \ref{as:terminal}. Assume for an induction hypothesis that  $Q_{t+1}^{i-1}(x) \geq V_{t+1}(x)$ for some $(i,t)\in I\setminus\{0\} \times T$ and for all $x\in X$. Fix any $x$ in $X$ and distinguish the two cases of the if-statement in step 12 of Algorithm \ref{alg:GBDP}.

\emph{Case I:}		
Suppose that $Q_{t+1}^{i-1}$ is submodular on $Z(x_{t+1}^i)$. Then $H^*$ is the unique hyperplane through the set $\{(y,(\mathcal{T}Q_{t+1}^{i})(y))\}_{y \in Y_+(x_{t+1}^i)}$. By \eqref{eq:qminh}, $Q_{t+1}^{i}$ is concave extensible since it is the pointwise minimum of a finite number of hyperplanes. Hence, we invoke Assumption \ref{as:preserve} to conclude that $\mathcal{T}Q_{t+1}^{i-1}$ is concave extensible and submodular. As shown by \cite[Appendix B.4]{LEBEDEVETAL2019B}, this implies that $H^*$ is a separating hyperplane, i.e.\ $H^*(x)\geq \mathcal{T}Q_{t+1}^{i-1}(x)$ for all $x\in X$. Define $d^V$ to be the maximiser of \eqref{eq:dp} and define $d^Q$ to be the maximiser of \eqref{eq:dp} with $V_{t+1}(y)$ replaced by $Q_{t+1}^{i-1}(y)$. We now show that the Bellman operator of the DP preserves the inequality $Q_{t+1}^{i-1}(x) \geq V_{t+1}(x)$, i.e.\ $\mathcal{T}Q_{t+1}^{i-1}(x) \geq \mathcal{T}V_{t+1}(x)$. To this end, fix $x\in X$ and consider				
\begin{align}
(\mathcal{T}Q_{t+1}^{i-1})(x)&= g(x,d^Q)+\sum_{y\in Y_+(x)}P_{x,y}(d^Q)Q_{t+1}^{i-1}(y)\nonumber\\
&\geq g(x,d^V)+\sum_{y\in Y_+(x)}P_{x,y}(d^V)Q_{t+1}^{i-1}(y)\nonumber\\
&\geq g(x,d^V)+\sum_{y\in Y_+(x)}P_{x,y}(d^V)V_{t+1}(y)\nonumber\\
&=(\mathcal{T}V_{t+1})(x),
\end{align}
where the first inequality follows from the supoptimality of $d^V$ for $(\mathcal{T}Q_{t+1}^{i-1})(x)$ and the second inequality follows from the induction hypothesis.
				
\emph{Case II:}
Now consider the case when $Q_{t+1}^{i-1}$ is not submodular on $Z(x_{t+1}^{i})$. Then $H^* \in \set{\mathcal{T}H_{t+1}^{j-1}}{j\in J_{t+1}^{i-1}}$. Furthermore, by \eqref{eq:qminh} and the induction hypothesis, 
\begin{equation}\label{eq:hqv}
H_{t+1}^{j-1}(x)\geq Q_{t+1}^{i-1}(x) \geq V_{t+1}(x), \text{ for all } (x,j)\in X\times J_{t+1}^{i-1}.
\end{equation}
We now show that all possible realisations of $H^*$ constitute upper bounds on $\mathcal{T}V_{t+1}$. To this end, fix any $(x,j)\in X\times J_{t+1}^{i-1}$. Define $d^H$ to be the maximiser of \eqref{eq:dp} with $V_{t+1}(y)$ replaced by $H_{t+1}^{j-1}(y)$. We can show that the Bellman operator of the DP preserves the inequality $Q_{t+1}^{i-1}(x) \geq V_{t+1}(x)$ using a similar argument as before:
\begin{equation}
(\mathcal{T}H_{t+1}^{j-1})(x)\geq(\mathcal{T}V_{t+1})(x),
\end{equation}
which follows from the suboptimality of $d^V$ (see Case I) for $(\mathcal{T}H_{t+1}^{j-1})(x)$ and the fact that $H_{t+1}^{j-1}(x)\geq V_{t+1}(x)$ (see \eqref{eq:hqv}). Therefore, we conclude that $H^*(x)\geq \mathcal{T}V_t(x)$ for all $x\in X$ in the second case as well. 
		
Since both cases lead to an upper bound, i.e.\ $H^*(x)\geq \mathcal{T}V_{t+1}(x)$ for all $x\in X$, we infer that 
\begin{equation}
Q_{t}^i(x)= \min\left\{H^*(x),Q_{t+1}^{i-1}(x)\right\} \geq  \mathcal{T}V_{t+1}(x)
\end{equation}
for all $x\in X$. This concludes our induction argument and shows that $Q_t^i(x) \geq V_t(x)$ for all $(x,i,t)\in X\times I \times T$.
\subsection{Proof of Proposition \ref{pr:converge}}
We will show the proposition by induction on $t$. Consider the base case, when $Q_{\bar{t}+1}^0(x) = V_{\bar{t}+1}(x)$ for all $x\in X$. Then notice that in the ``backward sweep'', the proposed algorithm computes the Bellman equation from $\bar{t}+1 \to \bar{t}$ exactly for every $x\in X$. This is because $Q_{\bar{t}+1}^0$ is submodular by Assumption \ref{as:terminal} and hence, the if-statement in step 12 of Algorithm \ref{alg:GBDP} is true. By Assumption \ref{as:resampling}, $x_{\bar{t}+1}^i$ is resampled if for the time step transition $\bar{t}+1 \to \bar{t}$, the algorithm has not converged to the exact value function at $\bar{t}$ yet. Therefore, the value function is computed exactly at all $x\in X$ for the time step transition $\bar{t}+1 \to \bar{t}$ after at most $|X|$ iterations of the proposed algorithm, i.e.\ $Q_{\bar{t}}^{\hat{i}}(x) = V_{\bar{t}}(x)$ for all $x\in X$, where $\hat{i}\leq |X|$.

Now suppose by means of an induction hypothesis that for some $(t,i)\in T\times I$,  $Q_{t+1}^i(x) = V_{t+1}(x)$ for all $x\in X$. Then by Assumptions \ref{as:terminal} and \ref{as:preserve}, $V_{t+1}$ is submodular and hence, $Q_{t+1}^i$ is also submodular. By a similar argument to the base case, the proposed algorithm computes the exact value function for the time step transition $t+1 \to t$ in another $\hat{i}\leq |X|$ iterations. 

Hence, we conclude that for every time step transition, the proposed algorithm needs at most $|X|$ iterations to compute the exact value function for any one time step $t\in T$, which gives at most $\bar{t}|X|$ iterations for the total time horizon. Hence, after any $i\geq\bar{t}|X|$ iterations, $Q_{t}^{i}(x) = V_{t}(x)$ for all $(x,t)\in X\times T$. Therefore, both $\mathbb{E}[l(i)]=V_1(0)$ and $u(i)=Q_{1}^{i}(0)=V_1(0)$, which finally implies that  $\mathbb{E}[l(i)]=u(i)$ for all $i\geq\bar{t}|X|$ iterations.

\subsection{Proof of Proposition \ref{pr:oneoffbounds}(i)}
The proof is a finite sample adaptation of the one-sided Chebyshev's inequality, i.e.\ Cantelli's inequality \cite[Theorem 1]{GHOSH2002}. We distinguish the following two cases:

\emph{Case I:} Suppose that $\sigma_{\mathrm{v}}\neq 0$. Fix any $k\in K$ and consider the conditional probability that $l_{\mathrm{C}}:=l_{\mathrm{v}}(k_{\max}+1)$ is no greater than $\bar{l}_{\mathrm{v}}-m\hat{\sigma}$ for some $m>0$:
\begin{equation}\label{eq:indicator}
\begin{split}
& \Pr(l_{\mathrm{C}}\leq \bar{l}_{\mathrm{v}}-m\sigma_{\mathrm{v}} | \sigma_{\mathrm{v}}\neq 0)\\
=& \Pr(m\sigma_{\mathrm{v}}\leq \bar{l}_{\mathrm{v}}-l_{\mathrm{C}}| \sigma_{\mathrm{v}}\neq 0)\\
=& \frac{1}{k_{\max}}\sum_{k\in K}\Pr(m\sigma_{\mathrm{v}}\leq \bar{l}_{\mathrm{v}}-l_{\mathrm{v}}(k)| \sigma_{\mathrm{v}}\neq 0)\\
=& \frac{1}{k_{\max}}\mathbb{E}\left(\left.\sum_{k\in K} \mathbbm{1}(m\sigma_{\mathrm{v}}\leq \bar{l}_{\mathrm{v}}-l_{\mathrm{v}}(k))\right|\sigma_{\mathrm{v}}\neq 0\right),
\end{split}
\end{equation}
where the second last equality follows from the observation that $l_{\mathrm{v}}(k)$ for all $k\in K\cup\{k_{\max}+1\}$ are independently and identically distributed. Next, we want to upper bound the indicator function in \eqref{eq:indicator} by a quadratic function. A suitable expression is given for any $c>0$ by
\begin{equation}
\begin{split}
& \Pr(l_{\mathrm{C}}\leq \bar{l}_{\mathrm{v}}-m\sigma_{\mathrm{v}} | \sigma_{\mathrm{v}}\neq 0)\\
\leq& \frac{1}{k_{\max}}\mathbb{E}\left(\left.\sum_{k\in K}\frac{(\bar{l}_{\mathrm{v}}-l_{\mathrm{v}}(k)+c\sigma_{\mathrm{v}})^2}{(m\sigma_{\mathrm{v}}+c\sigma_{\mathrm{v}})^2} \right| \sigma_{\mathrm{v}}\neq 0\right).
\end{split}
\end{equation}
Note that each element in the summation is always non-negative and no smaller than one if $m\sigma_{\mathrm{v}}\leq \bar{l}-l_{\mathrm{v}}(k)$ and hence is an upper bound to \eqref{eq:indicator}. We simplify this expression using the definitions of $\bar{l}_{\mathrm{v}}$ and $\sigma_{\mathrm{v}}$ from Algorithm \ref{alg:validate} in Section \ref{sec:valalg} as
\begin{equation}\label{eq:pr_bound_ito_c}
\begin{split}
& \frac{1}{k_{\max}}\mathbb{E}\left(\left.\sum_{k\in K}\frac{(\bar{l}_{\mathrm{v}}-l_{\mathrm{v}}(k)+c\sigma_{\mathrm{v}})^2}{(m\sigma_{\mathrm{v}}+c\sigma_{\mathrm{v}})^2}\right|\sigma_{\mathrm{v}}\neq 0\right)\\
=& \frac{1}{k_{\max}}\mathbb{E}\left(\left.\frac{(k_{\max}-1)\sigma_{\mathrm{v}}^2+k_{\max}c^2\sigma_{\mathrm{v}}^2}{(m+c)^2\sigma_{\mathrm{v}}^2}\right|\sigma_{\mathrm{v}}\neq 0\right)\\
=& \frac{1}{k_{\max}}\mathbb{E}\left(\frac{k_{\max}-1+k_{\max}c^2}{(m+c)^2}\right)\\
=& \frac{k_{\max}-1+k_{\max}c^2}{k_{\max}(m+c)^2},
\end{split}
\end{equation}
where $\sigma_{\mathrm{v}}$ cancels, since $\sigma_{\mathrm{v}}\neq 0$, and the expectation operator drops, since its argument is a constant. We minimise \eqref{eq:pr_bound_ito_c} by considering its first order condition, i.e.\
\begin{align}\label{eq:dbydc}
0=&\;\frac{\partial}{\partial c}\frac{k_{\max}-1+k_{\max}c^2}{k_{\max}(m+c)^2}\nonumber\\
=&\;(k_{\max}(m+c)^2)^{-2}\left\{2k_{\max}^2c(m+c)^2\right.\nonumber\\
&\;\left.-(k_{\max}-1+k_{\max}c^2)2k_{\max}(m+c)\right\}\nonumber\\
\Rightarrow c =&\; \frac{k_{\max}-1}{k_{\max}m}.
\end{align}
The second-order condition shows that $c$ minimises \eqref{eq:pr_bound_ito_c}. Substituting $c$ into \eqref{eq:pr_bound_ito_c} and simplifying gives:
\begin{equation}\label{eq:target_pr}
\Pr(l_{\mathrm{C}}\leq \bar{l}_{\mathrm{v}}-m\sigma_{\mathrm{v}} | \sigma_{\mathrm{v}}\neq 0)\leq\frac{k_{\max}-1}{k_{\max}m^2+k_{\max}-1}.
\end{equation}
\emph{Case II:} Suppose that $\sigma_{\mathrm{v}}=0$. We repeat the derivation of Case I with $\Pr(l_{\mathrm{C}}\leq \bar{l}_{\mathrm{v}}-m\sigma_{\mathrm{v}} | \sigma_{\mathrm{v}}\neq 0)$ replaced by $\Pr(l_{\mathrm{C}}\leq \bar{l}_{\mathrm{v}}-m\sigma_{\mathrm{v}} | \sigma_{\mathrm{v}}= 0)$ until \eqref{eq:indicator}, where we note that due to $\sigma_{\mathrm{v}}=0$, we have $l_{\mathrm{v}}(k)=l_{\mathrm{v}}(k')$ for all $(k,k')\in K\times K$ and hence, $\Pr(l_{\mathrm{C}}\leq \bar{l}_{\mathrm{v}}-m\sigma_{\mathrm{v}} | \sigma_{\mathrm{v}}=0)=1$.

Recalling that $\theta_{\mathrm{C}}:=\Pr(\sigma_{\mathrm{v}}=0)$ and taking both cases together, we obtain by the total probability theorem that
\begin{equation}
\Pr(l_{\mathrm{C}}\leq \bar{l}_{\mathrm{v}}-m\sigma_{\mathrm{v}})\leq\frac{(1-\theta_{\mathrm{C}})(k_{\max}-1)}{k_{\max}m^2+k_{\max}-1} + \theta_{\mathrm{C}}.
\end{equation}
We want this probability to be at most the significance level $\alpha$. Hence, we solve this expression for $m$, yielding
\begin{equation}
\begin{split}
\alpha&\geq(1-\theta_{\mathrm{C}}) \frac{k_{\max}-1}{k_{\max}m^2+k_{\max}-1} + \theta_{\mathrm{C}}\\
\iff m&\geq \sqrt{\frac{(1-\alpha)(k_{\max}-1)}{(\alpha -\theta_{\mathrm{C}})k_{\max}}},
\end{split}
\end{equation}
which is real-valued, since $\alpha > \theta_{\mathrm{C}}$ by Assumption \ref{as:alphac}. Substituting for $m$ in \eqref{eq:target_pr} gives the desired property:
\begin{equation}
\Pr\left(l_{\mathrm{C}}\leq \bar{l}-\sigma_{\mathrm{v}}\sqrt{\frac{(1-\alpha)(k_{\max}-1)}{(\alpha -\theta_{\mathrm{C}}) k_{\max}}}\right) \leq \alpha.
\end{equation}
\subsection{Proof of Proposition \ref{pr:oneoffbounds}(ii)}\label{sec:dkw_proof1}
Fix any $\alpha\in (0,1)$. By the Dvoretzky-Kiefer-Wolfowitz inequality \cite{MASSART1990}, we can write with (user-defined) probability at most $\theta_{\mathrm{D}}\in (0,\alpha)$ that
\begin{equation}\label{eq:dkw}
F(l_{\mathrm{D}})\geq F_K(l_{\mathrm{D}}) + \sqrt{\frac{\ln(1/\theta_{\mathrm{D}})}{2k_{\max}}},
\end{equation}
where $F$ is the exact yet unknown cumulative distribution function of $l_{\mathrm{v}}(k_{\max}+1)$, i.e.\ $\Pr(l_{\mathrm{v}}(k_{\max}+1)\leq l_{\mathrm{D}})$. Denote the random event in \eqref{eq:dkw} by $E$. Then, we can write
\begin{align}
\alpha=&\;\Pr(E \cup \Pr(l_{\mathrm{v}}(k_{\max}+1)\leq l_{\mathrm{D}})\nonumber \\
\leq &\; \Pr(E) + \Pr(l_{\mathrm{v}}(k_{\max}+1)\leq l_{\mathrm{D}})\nonumber \\
\leq&\; \theta_{\mathrm{D}}+F(l_{\mathrm{D}}),
\end{align}
where the first inequality is due to subadditivity of $\Pr$ (see \cite[Chapter 1.1]{HAJEK2015}). Hence, $F(l_{\mathrm{D}})\geq \alpha-\theta_{\mathrm{D}}$, which we substitute in \eqref{eq:dkw} to arrive at the required result
\begin{align}
& F_K(l_{\mathrm{D}}) \leq \alpha-\theta_{\mathrm{D}}- \sqrt{\frac{\ln(\frac{1}{\theta_{\mathrm{D}}})}{2k_{\max}}}\text{ and hence, }\nonumber\\
& l_{\mathrm{D}} = \mathrm{sup} \left\{l \in [l_-,l_+] \left| \vphantom{\sqrt{\frac{1}{1}}} F_K(l)\leq \alpha- \theta_{\mathrm{D}} \right.
-\sqrt{\frac{\ln(\frac{1}{\theta_{\mathrm{D}}})}{2k_{\max}}} \right\}.
\end{align}
\subsection{Proof of Lemma \ref{le:opt_ld}}
Consider the first-order optimality condition, i.e.\
\begin{align}
0&= \frac{\partial}{\partial \theta_{\mathrm{D}}}\left\{\alpha-\theta_{\mathrm{D}}-\sqrt{\frac{\ln(\frac{1}{\theta_{\mathrm{D}}})}{2k_{\max}}}\right\}\nonumber\\
\Rightarrow 0&= -1+\frac{1}{2\theta_{\mathrm{D}}\sqrt{-2k_{\max}\ln(\theta_{\mathrm{D}})}}.
\end{align}
Since $\theta_{\mathrm{D}}>0$, we can simplify to arrive at
\begin{align}
\theta_{\mathrm{D}}^2 \ln(\theta_{\mathrm{D}}^2) &= \frac{-1}{4k_{\max}}\nonumber\\
\Rightarrow \theta_{\mathrm{D}}^2 &= \exp\left(W_{i}\left(\frac{-1}{4k_{\max}}\right)\right),
\end{align}
where $i\in \{0,-1\}$ and $W_{i}$ is the Lambert $W$ function (see Definition \ref{de:lambertw}). The second-order conditions show that $i=-1$ gives a local maximum. Hence, we take the square root and note that $\theta_{\mathrm{D}}$ is bounded from above by $\alpha$ to arrive at the desired result.
\subsection{Proof of Proposition \ref{pr:ebounds}(ii)}
By the Dvoretzky-Kiefer-Wolfowitz inequality \cite{MASSART1990}, we can write with confidence $1-\alpha^{\mathbb{E}}$ that
\begin{equation}
F(l) \leq \min\left\{1,F_K(l) + \sqrt{\frac{\ln(1/\alpha^{\mathbb{E}})}{2k_{\max}}}\right\},
\end{equation}
where $F$ is defined as in Appendix \ref{sec:dkw_proof1}. Under Assumption \ref{as:lmin}, $l$ is non-negative. Hence, we find $\mathbb{E}l$ from $F$ as the following integral (see \cite[Chapter 1.5]{HAJEK2015}):
\begin{equation}
\begin{split}
1-F(l) \geq& 1-\min\left\{1,F_K(l) + \sqrt{\frac{\ln(1/\alpha^{\mathbb{E}})}{2k_{\max}}}\right\}\\
\Rightarrow \mathbb{E}l \geq & \int_{l=0}^{\infty}1-\min\left\{1,F_K(l) + \sqrt{\frac{\ln(1/\alpha^{\mathbb{E}})}{2k_{\max}}}\right\}\mathrm{d}l.
\end{split}
\end{equation}
Finally, $\mathbb{E}\bar{l}_{\mathrm{v}}=\mathbb{E}l$, since $\mathbb{E}$ is a linear operator and $l_{\mathrm{v}}(k)$ for all $k\in K$ are independent, thus concluding the proof.
\subsection{Proof of Proposition \ref{pr:DKW_H}}
Under Assumption \ref{as:lmin}, we can express the empirical mean in Hoeffding's bound \cite{HOEFFDING1963} as an integral over the empirical cumulative distribution function:
\begin{align}
l_{\mathrm{H}}^{\mathbb{E}}:=\;&\bar{l}_{\mathrm{v}}-(l_{+}-l_{-})\sqrt{\frac{\ln(1/\alpha^{\mathbb{E}})}{2k_{\max}}}\nonumber\\
=\;&\int_{l=0}^{\infty}1- F_K(l)\mathrm{d}l -(l_{+}-l_{-})\sqrt{\frac{\ln(1/\alpha^{\mathbb{E}})}{2k_{\max}}}\nonumber\\
\leq\;& \int_{l=l_{-}}^{l_{+}}1- \left(F_K(l)+\sqrt{\frac{\ln(1/\alpha^{\mathbb{E}})}{2k_{\max}}}\right)\mathrm{d}l\nonumber\\
<\;& \int_{l=l_{-}}^{l_{+}}1-\min\left\{1,F_K(l) + \sqrt{\frac{\ln(1/\alpha^{\mathbb{E}})}{2k_{\max}}}\right\}\mathrm{d}l\nonumber\\
=\;& l_{\mathrm{D}}^{\mathbb{E}},
\end{align}
where the last inequality is strict by Assumption \ref{as:alphae}, since $F_K$ is a stair function with step height $1/k_{\max}$ and $\sqrt{\ln(1/\alpha^{\mathbb{E}})/(2k_{\max})}>1/k_{\max}$ implies that $F_K(l^*)+\sqrt{\ln(1/\alpha^{\mathbb{E}})/(2k_{\max})} > 1$ for some $l^*<l_{+}$.
\end{document}